\documentclass[10pt]{article}
\usepackage{booktabs}
\usepackage{graphicx}
\usepackage{rotating}
\usepackage{mathtools}
\usepackage{latexsym}
\usepackage{amssymb}
\usepackage{mathrsfs}
\newtheorem{thm}{Theorem}[section]

\newtheorem{rmk}[thm]{Remark}

\newtheorem {defn}[thm]{Definition}

\begin{document}

\centerline {\large{\bf MADM Approach For Fermatean Neutrosophic Normal Aggregation Operator}}
\vspace{0.3cm}
\centerline{} 
 \centerline{\bf {M. Palanikumar$^1$, K. Arulmozhi$^2$, Santanu Acharjee$^{3, \dagger}$ }} 
\centerline{$^1$Department of Mathematics,}
\centerline{St.Joseph College of Engineering, Anna University, Sriperumbudur, Chennai-602117}
\centerline{$^2$Department of Mathematics,}
\centerline{Bharath Institute of Higher Education and Research, Chennai-600027}
\centerline{$^3$Department of Mathematics,}
\centerline{Gauhati University, Gauhati, Assam, India}
\centerline{e-mails: $^1$palanimaths86@gmail.com,$^2$arulmozhiems@gmail.com,}
\centerline{$^3$sacharjee326@gmail.com}
\centerline{$\dagger$corresponding author: Santanu Acharjee}
\centerline{$\dagger$Orcid id:  0000-0003-4932-3305}

\noindent {\bf Abstract:} 
We present a communication which deals with some new methods to solve multiple attribute decision-making (MADM) problems based on Fermatean neutrosophic normal  number (FNNN). Fermatean neutrosophic sets based on further generalization of neutrosophic and Pythagorean neutrosophic sets. To develop some Fermatean neutrosophic normal  aggregation operators. The notion of FNNN holds for commutative and associative laws. There are many aggregation operators that have been defined up to date, but we concentration of this article is to introduce a new concept of Fermatean neutrosophic normal  weighted averaging (FNNWA), Fermatean neutrosophic normal  weighted geometric(FNNWG), generalized Fermatean neutrosophic normal  weighted averaging(GFNNWA) and generalized Fermatean neutrosophic normal  weighted geometric(GFNNWG). Also, we obtain an algorithm that deals with the MADM problems based on these operators. We interact applicability of the euclidean and hamming distance measures which are further extended in real life example. Finally, some important properties of these sets under algebraic operations are to be elaborated in this communication.



\markboth{\centerline {\scriptsize M.Palanikumar, K.Arulmozhi and Santanu Acharjee}}
         {\centerline {\scriptsize MADM Approach For Fermatean Neutrosophic Normal Aggregation Operator Framework}}
\begin{section} {\bf Introduction} \label{intro}
\end{section}
\mbox{\hspace{1cm}} Uncertainty can be seen everywhere in most real problems. In order to cope with the uncertainties, many uncertain theories such as fuzzy set \cite{Zadeh}, intuitionistic fuzzy set \cite{Atanassov}, Pythagorean fuzzy set \cite{Yager2,Liang}, neutrosophic fuzzy set \cite{Smarandache1} and Fermatean fuzzy set \cite{Senapati1} are put forward. Zadeh was introduced by fuzzy set, which suggests that decision-makers are to be able to solve uncertain problems by considering membership degree. Later, the notion of an intuitionistic fuzzy set is launched by Atanassov \cite{Atanassov}. Yager \cite{Yager2} as having been introduced by the concept of Pythagorean fuzzy sets in 2014 and is classified by a degree of membership and non-membership that satisfies the condition that the square sum of its membership degree and non membership degree is not exceed unity. However, we may encounter a problem in decision-making events where the square sum of the degree of membership and non-membership of a particular attribute exceeds unity. In 2019, Senapati et al. proposed the notion of a Fermatean fuzzy set \cite{Senapati1}. It has been extended by the Pythagorean fuzzy sets and is characterised by the condition that the cubic sum of its degrees of membership and non-membership is not exceed unity. Rahman et al. \cite{Rahman1} worked on some geometric aggregation operators on interval-valued Pythagorean fuzzy sets and applied them to group decision-making problems. Rahman et al. \cite{Rahman2} proposed an approach to multi-attribute group decision making (MAGDM) based on an induced interval-valued Pythagorean fuzzy Einstein aggregation operator.\\ 
\mbox{\hspace{1cm}} The TOPSIS method consists of distances to positive ideal solutions (PIS) and negative ideal solutions (NIS), calculating a preference order that is ranked under relative closeness, and finding a combination of these two distance measures. In 2018, Akram and Arshad proposed the concept of bipolar fuzzy TOPSIS for group decision making (GDM) \cite{AkramArshad}. The $m$-polar fuzzy linguistic TOPSIS approach for multi-criteria group decision making (MCGDM) was discussed by Adeel et al. \cite{Adeel}. Zhang and Xu proposed the extension of Pythagorean fuzzy sets based on TOPSIS to multi- criteria decision making (MCDM) \cite{ZhangXu} and after 2018, Peng and Dai \cite{PengDai} initiated the concept of single-valued neutrosophic MADM based on MABAC and TOPSIS. Hwang et al. discussed various real life applications of MADM \cite{Hwang}. In 2019, Jana and Pal introduced the concept of Pythagorean fuzzy dombi aggregation operators \cite{Chiranjibe1}. In 2019, \cite{Ali} Ullah et al. discussed some distance measures of complex Pythagorean fuzzy sets and their applications in pattern recognition. In 2020, Jana et al.  trapezoidal neutrosophic aggregation operators and their application to the MADM process  \cite{Chiranjibe8}. In  2021, Jana and Pal interacted on the application for the MCDM process based on single valued neutrosophic dombi power aggregation operators \cite{Chiranjibe4}. In recent work, Jana et al. introduced the concept of an MCDM approach based on SVTrN dombi aggregation functions \cite{Chiranjibe7}. TOPSIS extends to interval-valued intuitionistic fuzzy soft set (IFSS) was discussed by Zulqarnain et al. in 2021. \\
\mbox{\hspace{1cm}} Recently, a new theory has been introduced which is known as neutrosophic logic and sets. The term "neutosophy" means knowledge of neutral thought, and this neutrality represents the main distinction between fuzzy and intuitionistic fuzzy logic and set. Neutrosophic logic was introduced by Florentin Smarandache \cite{Smarandache1}. It is a logic in which each proposition is estimated to have a degree of truth, a degree of indeterminacy  and a degree of falsity. A neutrosophic set is one in which each element of the universe has a degree of truth, indeterminacy, and falsity that ranges between $0$ and $1$. From the philosophical point of view, it has been shown that a neutrosophic set generalises a classical set, a fuzzy set, an interval valued fuzzy set, etc. Florentin Smarandache et al. introduced the Pythagorean neutrosophic interval-valued set (PNSIVS) \cite{Jansi}. An application of the single-valued neutrosophic set is based on medical diagnosis \cite{Shahzadi} and context analysis \cite{Singh}. Ejegwa \cite{Ejegwa} extended the distance measures for intuitionistic fuzzy sets, viz, hamming, euclidean, normalised hamming, and normalised euclidean distances, and similarities to Pythagorean fuzzy sets and applied them to MCDM problems as well as MADM problems. \\
\mbox{\hspace{1cm}} The purpose of this paper is to extend the concept of Fermatean neutrosophic (FN) information to Fermatean neutrosophic normal (FNN) information. We obtain FNN information from aggregation operators. The paper is organised into seven sections as follows. In Section \ref{intro} is called an introduction. In Section \ref{sec:2} brief description of FN information is given. In Section \ref{sec:3}, we discuss MADM based on FNNN and its operations. In Section \ref{sec:4}, we talk about MADM based on the modified distance method for FNNN. In Section \ref{sec:5}, we talk about MADM based on some aggregation operators for FNNN. In Section \ref{sec:6} MADM Based on FNN Information, an algorithm with a numerical example, analysis and discussion is presented. Finally, a conclusion is provided in Section \ref{sec:7}. Throughout this paper, we may interact $0 \leq (\tau^{\mathbb{T}}_{L}(u))^{2}+(\tau^{\mathbb{I}}_{L}(u))^{2} +(\tau^{\mathbb{F}}_{L}(u))^{2} > 2$ but $0 \leq (\tau^{\mathbb{T}}_{L}(u))^{3}+(\tau^{\mathbb{I}}_{L}(u))^{3} +(\tau^{\mathbb{F}}_{L}(u))^{3} \leq 2$.
Therefore, the main contributions of this work are:
\begin{enumerate}
\item We introduce an euclidean and a hamming distance for FNNNs.
\item Real-world numerical examples of the introduced concept's applicability for MADM and aggregation operators for FNNN
\item Determine the positive and negative ideal values for the concept for FNNWA, FNNWG, GFNNWA, and GFNNWG. 
\item Decision-making for arriving at $\Lambda$.
\end{enumerate}

\section{\bf Preliminaries}  \label{sec:2}
In this section, we discuss some existing definitions and results.

\begin{defn}  \cite{Yager2}
Let $\mathbb{U}$ be a universe of discourse. The PFS $L$ in $\mathbb{U}$ is given by $L= \Big\{u, \big\langle \tau^{\mathbb{T}}_{L}(u), \tau^{\mathbb{F}}_{L}(u)\big\rangle \big| u\in \mathbb{U}\Big\}$, where the functions $\tau^{\mathbb{T}}_{L}: \mathbb{U} \rightarrow [0,1]$ and $\tau^{\mathbb{F}}_{L}: \mathbb{U} \rightarrow [0,1]$  denote the membership degree and non-membership degree of the element $u\in \mathbb{U}$ to the set $L$, respectively and $0 \leq (\tau^{\mathbb{T}}_{L}(u))^{2}+(\tau^{\mathbb{F}}_{L}(u))^{2} \leq 1$. For convenience, $L = \big\langle  \tau^{\mathbb{T}}_{L},\tau_{L}^{\mathbb{F}} \big\rangle $ is called a Pythagorean fuzzy number(PFN).
\end{defn}

\begin{defn} \cite{Smarandache1}
Let $\mathbb{U}$ be a universe of discourse. The neutrosophic set(NS) $L$ in $\mathbb{U}$ is given by 
$L= \Big\{u, \big\langle \tau^{\mathbb{T}}_{L}(u), \tau^{\mathbb{I}}_{L}(u), \tau^{\mathbb{F}}_{L}(u)\big\rangle \big| u\in \mathbb{U}\Big\}$, where the functions $\tau^{\mathbb{T}}_{L}: \mathbb{U} \rightarrow [0,1]$, $\tau^{\mathbb{I}}_{L}: \mathbb{U} \rightarrow [0,1]$ and $\tau^{\mathbb{F}}_{L}: \mathbb{U} \rightarrow [0,1]$ denote the truth, indeterminacy and falsity membership degree of the element $u\in \mathbb{U}$ to the set $L$, respectively and  
$0 \leq \tau^{\mathbb{T}}_{L}(u) + \tau^{\mathbb{I}}_{L}(u)+ \tau^{\mathbb{F}}_{L}(u) \leq 3$. For convenience, $L= \big\langle  \tau^{\mathbb{T}}_{L},\tau^{\mathbb{I}}_{L},\tau_{L}^{\mathbb{F}} \big\rangle $ is called a neutrosophic number(NN).
\end{defn}

\begin{defn} \cite{Jansi}
Let $\mathbb{U}$ be a universe of discourse. The Pythagorean neutrosophic set(PNS) $L$ in $\mathbb{U}$ is given by $L= \Big\{u, \big\langle \tau^{\mathbb{T}}_{L}(u), \tau^{\mathbb{I}}_{L}(u), \tau^{\mathbb{F}}_{L}(u)\big\rangle \big| u\in \mathbb{U}\Big\}$, where the functions $\tau^{\mathbb{T}}_{L}: \mathbb{U} \rightarrow [0,1]$, $\tau^{\mathbb{I}}_{L}: \mathbb{U} \rightarrow [0,1]$ and $\tau^{\mathbb{F}}_{L}: \mathbb{U} \rightarrow [0,1]$ denote the truth, indeterminacy and falsity membership degree of the element $u\in \mathbb{U}$ to the set $L$, respectively and 
$0 \leq (\tau^{\mathbb{T}}_{L}(u))^{2}+(\tau^{\mathbb{I}}_{L}(u))^{2} +(\tau^{\mathbb{F}}_{L}(u))^{2} \leq 2$. For convenience, $ L= \big\langle  \tau^{\mathbb{T}}_{L},\tau^{\mathbb{I}}_{L},\tau_{L}^{\mathbb{F}} \big\rangle $ is called a Pythagorean neutrosophic number(PNN).
\end{defn}

\begin{defn}  \cite{Senapati1}
A Fermatean fuzzy set $L$ in the universe $\mathbb{U}$ is of the form : $L= \Big\{u, \big\langle \tau^{\mathbb{T}}_{L}(u), \tau^{\mathbb{F}}_{L}(u)\big\rangle \big| u\in \mathbb{U}\Big\}$, where $\tau^{\mathbb{T}}_{L}(u)$ and $\tau^{\mathbb{F}}_{L}(u)$ are called degree of membership  and degree of non membership of $u$ respectively. The mapping $\tau^{\mathbb{T}}_{L}: \mathbb{U} \rightarrow [0,1]$ and  $\tau^{\mathbb{F}}_{L} : \mathbb{U} \rightarrow [0,1]$ and $0 \leq (\tau^{\mathbb{T}}_{L}(u))^{3}+(\tau^{\mathbb{F}}_{L}(u))^{3} \leq 1$. Since $L= \big\langle  \tau^{\mathbb{T}}_{L},\tau_{L}^{\mathbb{F}} \big\rangle $ is called a Fermatean fuzzy number(FFN). 
\end{defn}

\begin{defn} \cite{Senapati1} \label{defmain}
For any FFNs, $L= \langle  \tau^{\mathbb{T}} ,\tau^{\mathbb{F}}\rangle $, $L_{1}= \langle  \tau^{\mathbb{T}}_{1}, \tau^{\mathbb{F}}_{1} \rangle $ and $L_{2}= \langle  \tau^{\mathbb{T}}_{2}, \tau^{\mathbb{F}}_{2} \rangle $, $\tau^{\mathbb{T}},\tau^{\mathbb{F}}$ are called degree of membership and degree of non membership of $u$ respectively and $\Lambda$ be a positive integer. Then, the following statements are valid.
\begin{enumerate}
\item 
$L_{1}\boxplus L_{2}=
\begin{pmatrix}
\sqrt[3]{(\tau^{\mathbb{T}}_{1})^{3} + (\tau^{\mathbb{T}}_{2})^{3} -(\tau^{\mathbb{T}}_{1})^{3} \cdot (\tau^{\mathbb{T}}_{2})^{3}},
(\tau^{\mathbb{F}}_{1}\cdot\tau^{\mathbb{F}}_{2})
\end{pmatrix}$
\item
$L_{1} \boxtimes L_{2} =  \begin{pmatrix}
(\tau^{\mathbb{T}}_{1}\cdot \tau^{\mathbb{T}}_{2}),
\sqrt[3]{(\tau^{\mathbb{F}}_{1})^{3} + (\tau^{\mathbb{F}}_{2})^{3} -(\tau^{\mathbb{F}}_{1})^{3} \cdot (\tau^{\mathbb{F}}_{2})^{3}}
\end{pmatrix}$
\item $ \Lambda \cdot L = \begin{pmatrix}
\sqrt[3]{1-\big( 1- (\tau^{\mathbb{T}})^{3}\big)^{\Lambda}},
(\tau^{\mathbb{F}})^{\Lambda}
\end{pmatrix}$
\item $L^{\Lambda}= \begin{pmatrix}
(\tau^{\mathbb{T}})^{\Lambda},
\sqrt[3]{1-\big( 1- (\tau^{\mathbb{F}})^{3}\big)^{\Lambda}}\,\,
\end{pmatrix}$
\end{enumerate} 
\end{defn}

\begin{defn} \cite{Senapati1}
Let $L= \big\langle  \tau^{\mathbb{T}},  \tau^{\mathbb{F}}  \big\rangle $ be FFN, then its score function is defined as
$$S(L)= (\tau^{\mathbb{T}})^{3} - (\tau^{\mathbb{F}})^{3},\,\, S(L)\in [-1,1]$$
and its accuracy function is defined as
$$H(L)= (\tau^{\mathbb{T}})^{3}   +(\tau^{\mathbb{F}})^{3},\,\, H(L)\in [0,1].$$
\end{defn} 

\begin{defn} \cite{Yang1}
Let $R$ be a real number set, the membership function of fuzzy number $M(x)= e^{-\big(\frac{x-\eta}{\xi}\big)^{2}},(\xi >0)$ is called a normal fuzzy number (NFN) $M= (\eta, \xi)$ and the normal fuzzy number (NFN) set is denoted by $\mathbb{N}$.
\end{defn}

\section{Operations for FNNN} \label{sec:3}
In this section, we characterize the concept of FNN with Fermatean normal fuzzy number.

\begin{defn} 
A neutrosophic fuzzy set $L$ in the universe $\mathbb{U}$ is of the form $L= \Big\{u, \big\langle \tau^{\mathbb{T}}_{L}(u), \tau^{\mathbb{I}}_{L}(u), \tau^{\mathbb{F}}_{L}(u)\big\rangle \big| u\in \mathbb{U}\Big\}$, where $\tau^{\mathbb{T}}_{L}(u)$, $\tau^{\mathbb{I}}_{L}(u)$ and $\tau^{\mathbb{F}}_{L}(u)$ are called degree of truth, indeterminacy and falsity-membership of $L$ respectively. The mapping $\tau^{\mathbb{T}}_{L}: \mathbb{U} \rightarrow [0,1]$, $\tau^{\mathbb{I}}_{L} : \mathbb{U} \rightarrow [0,1]$, $\tau^{\mathbb{F}}_{L} : \mathbb{U} \rightarrow [0,1]$ and $0 \leq (\tau^{\mathbb{T}}_{L}(u))^{3}+(\tau^{\mathbb{I}}_{L}(u))^{3} +(\tau^{\mathbb{F}}_{L}(u))^{3} \leq 2$. Here, $L= \big\langle  \tau^{\mathbb{T}}_{L},\tau^{\mathbb{I}}_{L},\tau_{L}^{\mathbb{F}} \big\rangle $ is called a Fermatean neutrosophic number (FNN). 
\end{defn}

\begin{defn} 
\begin{enumerate}
\item Let $R$ be the real numbers, the Fermatean fuzzy number $Q(x)= e^{-\big(\frac{x-\eta}{\xi}\big)^{3}},(\xi >0)$ is called a Fermatean normal fuzzy number $Q= (\eta, \xi)$. The set of Fermatean normal fuzzy numbers is denoted by $\mathbb{N}$.
\item Let $P=(\eta_{1},\xi_{1})\in \mathbb{N}$and $Q=(\eta_{2},\xi_{2})\in \mathbb{N}$, then the distance between $P$ and $Q$ can be defined as $\mathbb{D}(P,Q)= \sqrt[3]{(\eta_{1}-\eta_{2})^{3} + \frac{1}{3}(\xi_{1}-\xi_{2})^{3}}$.
\end{enumerate}
\end{defn}

\begin{defn} 
Let $(\eta, \xi) \in \mathbb{N}$, $L= \big\langle (\eta, \xi); \tau^{\mathbb{T}},\tau^{\mathbb{I}}, \tau^{\mathbb{F}}\big\rangle  $ is a Fermatean neutrosophic normal number (FNNN), where its degree of truth, indeterminacy and falsity membership are defined as $\tau_{L}^{\mathbb{T}}= \tau_{L}^{\mathbb{T}} e^{-\big(\frac{x-\eta}{\xi}\big)^{3}}$, $\tau_{L}^{\mathbb{I}}= \tau_{L}^{\mathbb{I}} e^{-\big(\frac{x-\eta}{\xi}\big)^{3}}$ and $\tau_{L}^{\mathbb{F}}= 1-\big(1-\tau_{L}^{\mathbb{F}}\big) e^{-\big(\frac{x-\eta}{\xi}\big)^{3}}$, $x \in X$ respectively, where $X$ is a non-empty set and $\tau_{L}^{\mathbb{T}},\tau_{L}^{\mathbb{I}}, \tau_{L}^{\mathbb{F}}\in [0,1]$ and $0 \leq \big(\tau_{L}^{\mathbb{T}}(x)\big)^{3} + \big(\tau_{L}^{\mathbb{I}}(x)\big)^{3}+ \big(\tau_{L}^{\mathbb{F}}(x)\big)^{3} \leq 2$. 
\end{defn} 

\begin{defn} \label{def2}
If $L_{1}= \big\langle  (\eta_{1}, \xi_{1});  \tau^{\mathbb{T}}_{1},\tau^{\mathbb{I}}_{1} , \tau^{\mathbb{F}}_{1}  \big\rangle $, $L_{2}= \big\langle  (\eta_{2}, \xi_{2});  \tau^{\mathbb{T}}_{2}, \tau^{\mathbb{I}}_{2},\tau^{\mathbb{F}}_{2}  \big\rangle $ be the any two FNNNs and  $\Lambda$ be a positive integer, then
\begin{enumerate}
\item 
$L_{1}\boxplus L_{2} =
\begin{pmatrix}
(\eta_{1} + \eta_{2}, \xi_{1} + \xi_{2}); \sqrt[3\Lambda]{(\tau^{\mathbb{T}}_{1})^{3\Lambda} + (\tau^{\mathbb{T}}_{2})^{3\Lambda} -(\tau^{\mathbb{T}}_{1})^{3\Lambda} \cdot (\tau^{\mathbb{T}}_{2})^{3\Lambda}},\\
\sqrt[\Lambda]{(\tau^{\mathbb{I}}_{1})^{\Lambda} + (\tau^{\mathbb{I}}_{2})^{\Lambda} -(\tau^{\mathbb{I}}_{1})^{\Lambda} \cdot (\tau^{\mathbb{I}}_{2})^{\Lambda}},(\tau^{\mathbb{F}}_{1}\cdot \tau^{\mathbb{F}}_{2})
\end{pmatrix}$
\item 
$L_{1} \boxtimes L_{2} =  
\begin{pmatrix}
(\eta_{1} \cdot \eta_{2}, \xi_{1} \cdot \xi_{2});
(\tau^{\mathbb{T}}_{1}\cdot \tau^{\mathbb{T}}_{2}),
\sqrt[\Lambda]{(\tau^{\mathbb{I}}_{1})^{\Lambda} + (\tau^{\mathbb{I}}_{2})^{\Lambda} -(\tau^{\mathbb{I}}_{1})^{\Lambda} \cdot (\tau^{\mathbb{I}}_{2})^{\Lambda}},\\
\sqrt[3\Lambda]{(\tau^{\mathbb{F}}_{1})^{3\Lambda} + (\tau^{\mathbb{F}}_{2})^{3\Lambda} -(\tau^{\mathbb{F}}_{1})^{3\Lambda} \cdot (\tau^{\mathbb{F}}_{2})^{3\Lambda}}
\end{pmatrix}$
\item $ \Lambda \cdot L_{1} = \begin{pmatrix}
(\Lambda \cdot \eta_{1}, \Lambda \cdot \xi_{1}); \sqrt[3\Lambda]{1-\big( 1- (\tau^{\mathbb{T}}_{1})^{3\Lambda}\big)^{\Lambda}}, (\tau^{\mathbb{I}}_{1})^{\Lambda}, (\tau^{\mathbb{F}}_{1})^{\Lambda}
\end{pmatrix}$
\item $ L_{1}^{\Lambda} = \begin{pmatrix}
(\eta_{1}^{\Lambda}, \xi_{1}^{\Lambda}); (\tau^{\mathbb{T}}_{1})^{\Lambda},
(\tau^{\mathbb{I}}_{1})^{\Lambda},
\sqrt[3\Lambda]{1-\big( 1- (\tau^{\mathbb{F}}_{1})^{3\Lambda}\big)^{\Lambda}}\,\,
\end{pmatrix}$
\end{enumerate}
\end{defn} 

Now, we characterize the following theorem under the above operator.

\begin{thm}
Let $L_{1}= \big\langle  (\eta_{1}, \xi_{1});  \tau^{\mathbb{T}}_{1},\tau^{\mathbb{I}}_{1}, \tau^{\mathbb{F}}_{1}  \big\rangle $, $L_{2}= \big\langle  (\eta_{2}, \xi_{2});  \tau^{\mathbb{T}}_{2},\tau^{\mathbb{I}}_{2}, \tau^{\mathbb{F}}_{2}  \big\rangle $ and $L_{3}= \big\langle  (\eta_{3}, \xi_{3});  \tau^{\mathbb{T}}_{3}, \tau^{\mathbb{I}}_{3},\tau^{\mathbb{F}}_{3}  \big\rangle $ be any three FNNNs. Then the following statements are valid.
\begin{enumerate}
\item $L_{1} \boxplus L_{2} = L_{2} \boxplus L_{1}$
\item $(L_{1} \boxplus L_{2})\boxplus {L_{3}} = L_{1} \boxplus (L_{2}\boxplus {L_{3}})$
\item $L_{1} \boxtimes L_{2} = L_{2} \boxtimes L_{1}$
\item $(L_{1} \boxtimes L_{2})\boxtimes {L_{3}} = L_{1} \boxtimes (L_{2}\boxtimes {L_{3}})$
\end{enumerate}
\end{thm} 
\textbf{Proof.} The proofs of (1) and (3) are straightforward.\\
Now, we prove (2), Let $\mathbb{N}_{L}$ be the normal fuzzy number of FNN $L$, the degree of truth membership of $(L_{1} \boxplus L_{2})\boxplus {L_{3}}$ and $L_{1} \boxplus (L_{2}\boxplus {L_{3}})$ be $\tau^{\mathbb{T}}_{(L_{1} \boxplus L_{2})\boxplus L_{3}}$ and 
$\tau^{\mathbb{T}}_{L_{1} \boxplus (L_{2}\boxplus L_{3})}$ respectively, the degree of indeterminacy membership be $\tau^{\mathbb{I}}_{(L_{1} \boxplus L_{2})\boxplus L_{3}}$ and $\tau^{\mathbb{I}}_{L_{1} \boxplus (L_{2}\boxplus L_{3})}$ respectively, and the degree of falsity membership  be $\tau^{\mathbb{F}}_{(L_{1} \boxplus L_{2})\boxplus L_{3}}$ and $\tau^{\mathbb{F}}_{L_{1} \boxplus (L_{2}\boxplus L_{3})}$ respectively.\\
We get, $\mathbb{N}_{(L_{1} \boxplus L_{2})\boxplus L_{3}}= \mathbb{N}_{L_{1} \boxplus (L_{2}\boxplus L_{3})}= (\eta_{1}+\eta_{2}+\eta_{3}, \xi_{1}+\xi_{2}+\xi_{3})$.\\
Now,
\begin{eqnarray} \label{3.1}
\tau^{\mathbb{T}}_{(L_{1} \boxplus L_{2})\boxplus L_{3}}
&=& \nonumber
\sqrt[3\Lambda]{
\begin{aligned}
\Big((\tau^{\mathbb{T}}_{1})^{3\Lambda} + (\tau^{\mathbb{T}}_{2})^{3\Lambda} -(\tau^{\mathbb{T}}_{1})^{3\Lambda} \cdot (\tau^{\mathbb{T}}_{2})^{3\Lambda}\Big) + (\tau^{\mathbb{T}}_{3})^{3\Lambda}\\
-\Big((\tau^{\mathbb{T}}_{1})^{3\Lambda} + (\tau^{\mathbb{T}}_{2})^{3\Lambda} -(\tau^{\mathbb{T}}_{1})^{3\Lambda} \cdot (\tau^{\mathbb{T}}_{2})^{3\Lambda}\Big) \cdot (\tau^{\mathbb{T}}_{3})^{3\Lambda}\\  
\end{aligned}}\\
&=&
\sqrt[3\Lambda]{
\begin{aligned}
(\tau^{\mathbb{T}}_{1})^{3\Lambda} + (\tau^{\mathbb{T}}_{2})^{3\Lambda}+ (\tau^{\mathbb{T}}_{3})^{3\Lambda} -(\tau^{\mathbb{T}}_{1})^{3\Lambda} \cdot (\tau^{\mathbb{T}}_{2})^{3\Lambda} -(\tau^{\mathbb{T}}_{1})^{3\Lambda} \cdot (\tau^{\mathbb{T}}_{3})^{3\Lambda}\\
-(\tau^{\mathbb{T}}_{2})^{3\Lambda} \cdot (\tau^{\mathbb{T}}_{3})^{3\Lambda} +(\tau^{\mathbb{T}}_{1})^{3\Lambda} \cdot (\tau^{\mathbb{T}}_{2})^{3\Lambda}\cdot (\tau^{\mathbb{T}}_{3})^{3\Lambda}
\end{aligned}}
\end{eqnarray}
and 
\begin{eqnarray} \label{3.2}
\tau^{\mathbb{T}}_{L_{1} \boxplus (L_{2}\boxplus L_{3})}
&=& \nonumber
\sqrt[3\Lambda]{
\begin{aligned}
\Big((\tau^{\mathbb{T}}_{2})^{3\Lambda} + (\tau^{\mathbb{T}}_{3})^{3\Lambda} -(\tau^{\mathbb{T}}_{2})^{3\Lambda} \cdot (\tau^{\mathbb{T}}_{3})^{3\Lambda}\Big) + (\tau^{\mathbb{T}}_{1})^{3\Lambda}\\
-\Big((\tau^{\mathbb{T}}_{2})^{3\Lambda} + (\tau^{\mathbb{T}}_{3})^{3\Lambda} -(\tau^{\mathbb{T}}_{2})^{3\Lambda} \cdot (\tau^{\mathbb{T}}_{3})^{3\Lambda}\Big) \cdot (\tau^{\mathbb{T}}_{1})^{3\Lambda}
\end{aligned}} \\
&=& 
\sqrt[3\Lambda]{
\begin{aligned}
(\tau^{\mathbb{T}}_{1})^{3\Lambda} + (\tau^{\mathbb{T}}_{2})^{3\Lambda}+ (\tau^{\mathbb{T}}_{3})^{3\Lambda} -(\tau^{\mathbb{T}}_{1})^{3\Lambda} \cdot (\tau^{\mathbb{T}}_{2})^{3\Lambda} -(\tau^{\mathbb{T}}_{1})^{3\Lambda} \cdot (\tau^{\mathbb{T}}_{3})^{3\Lambda}\\
 -(\tau^{\mathbb{T}}_{2})^{3\Lambda} \cdot (\tau^{\mathbb{T}}_{3})^{3\Lambda} +(\tau^{\mathbb{T}}_{1})^{3\Lambda} \cdot (\tau^{\mathbb{T}}_{2})^{3\Lambda}\cdot (\tau^{\mathbb{T}}_{3})^{3\Lambda}
\end{aligned}}
\end{eqnarray}
From the equations \ref{3.1} and Equation \ref{3.2}, $\tau^{\mathbb{T}}_{(L_{1} \boxplus L_{2})\boxplus L_{3}} = \tau^{\mathbb{T}}_{L_{1} \boxplus (L_{2}\boxplus L_{3})}$. 
Similarly, we  can prove that $\tau^{\mathbb{I}}_{(L_{1} \boxplus L_{2})\boxplus L_{3}} = \tau^{\mathbb{I}}_{L_{1} \boxplus (L_{2}\boxplus L_{3})}$ and  $\tau^{\mathbb{F}}_{(L_{1} \boxplus L_{2})\boxplus L_{3}} = \tau^{\mathbb{F}}_{L_{1} \boxplus (L_{2}\boxplus L_{3})}$.
Hence, $(L_{1} \boxplus L_{2})\boxplus {L_{3}} = L_{1} \boxplus (L_{2}\boxplus {L_{3}})$.\\
Let the degree of truth, degree of indeterminacy, and degree of falsity  membership of $(L_{1} \boxtimes L_{2})\boxtimes {L_{3}}$ and $L_{1} \boxtimes (L_{2}\boxtimes {L_{3}})$ be $\tau^{\mathbb{T}}_{(L_{1} \boxtimes L_{2})\boxtimes L_{3}} $ and $\tau^{\mathbb{T}}_{L_{1} \boxtimes (L_{2}\boxtimes L_{3})}$, $\tau^{\mathbb{I}}_{(L_{1} \boxtimes L_{2})\boxtimes L_{3}}$ and $\tau^{\mathbb{I}}_{L_{1} \boxtimes (L_{2}\boxtimes L_{3})}$, $\tau^{\mathbb{F}}_{(L_{1} \boxtimes L_{2})\boxtimes L_{3}}$ and $\tau^{\mathbb{F}}_{L_{1} \boxtimes (L_{2}\boxtimes L_{3})}$ respectively. We get, $\mathbb{N}_{(L_{1} \boxtimes L_{2})\boxtimes L_{3}}= \mathbb{N}_{L_{1} \boxtimes (L_{2}\boxtimes L_{3})}= (\eta_{1} \cdot \eta_{2} \cdot \eta_{3}, \xi_{1} \cdot \xi_{2} \cdot \xi_{3})$.\\
Now,
\begin{eqnarray} \label{3.3}
\tau^{\mathbb{F}}_{(L_{1} \boxtimes L_{2})\boxtimes L_{3}} 
&=& \nonumber
\sqrt[3\Lambda]{
\begin{aligned}
\Big((\tau^{\mathbb{F}}_{1})^{3\Lambda} + (\tau^{\mathbb{F}}_{2})^{3\Lambda} -(\tau^{\mathbb{F}}_{1})^{3\Lambda} \cdot (\tau^{\mathbb{F}}_{2})^{3\Lambda}\Big) + (\tau^{\mathbb{F}}_{3})^{3\Lambda}\\
-\Big((\tau^{\mathbb{F}}_{1})^{3\Lambda} + (\tau^{\mathbb{F}}_{2})^{3\Lambda} -(\tau^{\mathbb{F}}_{1})^{3\Lambda} \cdot (\tau^{\mathbb{F}}_{2})^{3\Lambda}\Big) \cdot (\tau^{\mathbb{F}}_{3})^{3\Lambda}
\end{aligned}}\\
&=& \sqrt[3\Lambda]{
\begin{aligned}
(\tau^{\mathbb{F}}_{1})^{3\Lambda} + (\tau^{\mathbb{F}}_{2})^{3\Lambda}+ (\tau^{\mathbb{F}}_{3})^{3\Lambda} -(\tau^{\mathbb{F}}_{1})^{3\Lambda} \cdot (\tau^{\mathbb{F}}_{2})^{3\Lambda} -(\tau^{\mathbb{F}}_{1})^{3\Lambda} \cdot (\tau^{\mathbb{F}}_{3})^{3\Lambda}\\
-(\tau^{\mathbb{F}}_{2})^{3\Lambda} \cdot (\tau^{\mathbb{F}}_{3})^{3\Lambda} +(\tau^{\mathbb{F}}_{1})^{3\Lambda} \cdot (\tau^{\mathbb{F}}_{2})^{3\Lambda}\cdot (\tau^{\mathbb{F}}_{3})^{3\Lambda}
\end{aligned}}
\end{eqnarray} 
\begin{eqnarray} \label{3.4}
\tau^{\mathbb{F}}_{L_{1} \boxtimes (L_{2}\boxtimes L_{3})}
&=& \nonumber
\sqrt[3\Lambda]{
\begin{aligned}
\Big((\tau^{\mathbb{F}}_{2})^{3\Lambda} + (\tau^{\mathbb{F}}_{3})^{3\Lambda} -(\tau^{\mathbb{F}}_{2})^{3\Lambda} \cdot (\tau^{\mathbb{F}}_{3})^{3\Lambda}\Big) + (\tau^{\mathbb{F}}_{1})^{3\Lambda}\\
-\Big((\tau^{\mathbb{F}}_{2})^{3\Lambda} + (\tau^{\mathbb{F}}_{3})^{3\Lambda} -(\tau^{\mathbb{F}}_{2})^{3\Lambda} \cdot (\tau^{\mathbb{F}}_{3})^{3\Lambda}\Big) \cdot (\tau^{\mathbb{F}}_{1})^{3\Lambda}
\end{aligned}}\\
&=& \sqrt[3\Lambda]{
\begin{aligned}
(\tau^{\mathbb{F}}_{1})^{3\Lambda} + (\tau^{\mathbb{F}}_{2})^{3\Lambda}+ (\tau^{\mathbb{F}}_{3})^{3\Lambda} -(\tau^{\mathbb{F}}_{1})^{3\Lambda} \cdot (\tau^{\mathbb{F}}_{2})^{3\Lambda} -(\tau^{\mathbb{F}}_{1})^{3\Lambda} \cdot (\tau^{\mathbb{F}}_{3})^{3\Lambda}\\
 -(\tau^{\mathbb{F}}_{2})^{3\Lambda} \cdot (\tau^{\mathbb{F}}_{3})^{3\Lambda} +(\tau^{\mathbb{F}}_{1})^{3\Lambda} \cdot (\tau^{\mathbb{F}}_{2})^{3\Lambda}\cdot (\tau^{\mathbb{F}}_{3})^{3\Lambda}
\end{aligned}}
\end{eqnarray}
From equations \ref{3.3} and \ref{3.4},  we have, $\tau^{\mathbb{F}}_{(L_{1} \boxtimes L_{2})\boxtimes L_{3}} = \tau^{\mathbb{F}}_{L_{1} \boxtimes (L_{2}\boxtimes L_{3})}$. Similarly, we can prove that $\tau^{\mathbb{T}}_{(L_{1} \boxtimes L_{2})\boxtimes L_{3}} = \tau^{\mathbb{T}}_{L_{1} \boxtimes (L_{2}\boxtimes L_{3})}$  and $\tau^{\mathbb{I}}_{(L_{1} \boxtimes L_{2})\boxtimes L_{3}} = \tau^{\mathbb{I}}_{L_{1} \boxtimes (L_{2}\boxtimes L_{3})}$. Hence $(L_{1} \boxtimes L_{2})\boxtimes {L_{3}} = L_{1} \boxtimes (L_{2}\boxtimes {L_{3}})$.

\section{Distance measure for FNNN} \label{sec:4}
In this section, we discuss about distance measure between two FNNNs. Distance measures plays a new crucial role in soft compiling to find the distances between the two methods lets with fuzzy, IFS,  neutrosophic, etc. 

\begin{defn} \label{def4.1}
Let $L_{1}= \big\langle  (\eta_{1}, \xi_{1});  \tau^{\mathbb{T}}_{1},\tau^{\mathbb{I}}_{1}, \tau^{\mathbb{F}}_{1}  \big\rangle $ and  $L_{2}= \big\langle  (\eta_{2}, \xi_{2});  \tau^{\mathbb{T}}_{2},\tau^{\mathbb{I}}_{2}, \tau^{\mathbb{F}}_{2}  \big\rangle $ be the any two FNNNs. Then, \\
(i) the Euclidean distance between $L_{1}$ and $L_{2}$ is defined as $\mathbb{D}_{E}(L_{1}, L_{2})$, where
$$\mathbb{D}_{E}(L_{1}, L_{2})=
\frac{1}{3}
\sqrt[3]{
\begin{aligned}
\begin{pmatrix}
\frac{1+(\tau^{\mathbb{T}}_{1})^{3}+(\tau^{\mathbb{I}}_{1})^{3}- (\tau^{\mathbb{F}}_{1})^{3}}{3} \eta_{1}
- \frac{1+(\tau^{\mathbb{T}}_{2})^{3}+(\tau^{\mathbb{I}}_{2})^{3}-(\tau^{\mathbb{F}}_{2})^{3}}{3} \eta_{2}
\end{pmatrix}^{3}\\
+ \frac{1}{3} 
\begin{pmatrix}
\frac{1+(\tau^{\mathbb{T}}_{1})^{3}+(\tau^{\mathbb{I}}_{1})^{3}-(\tau^{\mathbb{F}}_{1})^{3} }{3} \xi_{1}
- \frac{1+(\tau^{\mathbb{T}}_{2})^{3}+(\tau^{\mathbb{I}}_{2})^{3}-(\tau^{\mathbb{F}}_{2})^{3}}{3} \xi_{2} 
\end{pmatrix}^{3}
\end{aligned}}$$
(ii) the Hamming distance between $L_{1}$ and $L_{2}$ is defined as $\mathbb{D}_{H}(L_{1}, L_{2})$, where
$$\mathbb{D}_{H}(L_{1}, L_{2})=
\frac{1}{3}
\begin{pmatrix}
\left|
\begin{aligned}
\frac{1+(\tau^{\mathbb{T}}_{1})^{3}+(\tau^{\mathbb{I}}_{1})^{3}- (\tau^{\mathbb{F}}_{1})^{3}}{3} \eta_{1}
- \frac{1+(\tau^{\mathbb{T}}_{2})^{3}+(\tau^{\mathbb{I}}_{2})^{3}-(\tau^{\mathbb{F}}_{2})^{3} }{3} \eta_{2}
\end{aligned}
\right|\\
+ \frac{1}{3} 
\left|
\begin{aligned}
\frac{1+(\tau^{\mathbb{T}}_{1})^{3}+(\tau^{\mathbb{I}}_{1})^{3}-(\tau^{\mathbb{F}}_{1})^{3}}{3} \xi_{1}
- \frac{1+(\tau^{\mathbb{T}}_{2})^{3}+(\tau^{\mathbb{I}}_{2})^{3}- (\tau^{\mathbb{F}}_{2})^{3}}{3} \xi_{2}
\end{aligned}
\right|\\
\end{pmatrix}$$
\end{defn}

\begin{thm} \label{ed1}
Let $L_{1} = \big\langle  (\eta_{1}, \xi_{1}); \tau^{\mathbb{T}}_{1},\tau^{\mathbb{I}}_{1}, \tau^{\mathbb{F}}_{1}  \big\rangle$, $L_{2} = \big\langle  (\eta_{2}, \xi_{2}); \tau^{\mathbb{T}}_{2}, \tau^{\mathbb{I}}_{2},\tau^{\mathbb{F}}_{2}  \big\rangle  $ and ${L_{3}} = \big\langle  (\eta_{3}, \xi_{3}); \tau^{\mathbb{T}}_{3}, \tau^{\mathbb{I}}_{3}, \tau^{\mathbb{F}}_{3}  \big\rangle $ be any three FNNNs, then the following properties are hold.
\begin{enumerate}
\item $\mathbb{D}_{E}(L_{1}, L_{2}) = 0$ if and only if $L_{1}= L_{2}$
\item $ \mathbb{D}_{E}(L_{1}, {L_{3}}) \leq \mathbb{D}_{E}(L_{1}, L_{2}) + \mathbb{D}_{E}(L_{2}, {L_{3}})$.
\end{enumerate}
\end{thm}
\textbf{Proof.} The proof of (1) is straightforward. We will prove (2). Now,$(\mathbb{D}_{E}(L_{1}, L_{2}) + \mathbb{D}_{E}(L_{2}, L_{3}))^{3}$ can be written as
{\small
\begin{eqnarray*}
&& \frac{1}{27}
\begin{pmatrix}
\sqrt[3]{
\begin{aligned}
\begin{pmatrix}
\frac{1+(\tau^{\mathbb{T}}_{1})^{3}+(\tau^{\mathbb{I}}_{1})^{3}- (\tau^{\mathbb{F}}_{1})^{3}}{3} \eta_{1}
- \frac{1+(\tau^{\mathbb{T}}_{2})^{3}+(\tau^{\mathbb{I}}_{2})^{3}-(\tau^{\mathbb{F}}_{2})^{3} }{3} \eta_{2}
\end{pmatrix}^{3}\\
+ \frac{1}{3} 
\begin{pmatrix}
\frac{1+(\tau^{\mathbb{T}}_{1})^{3}+(\tau^{\mathbb{I}}_{1})^{3}-(\tau^{\mathbb{F}}_{1})^{3} }{3} \xi_{1}
- \frac{1+(\tau^{\mathbb{T}}_{2})^{3}+(\tau^{\mathbb{I}}_{2})^{3}- (\tau^{\mathbb{F}}_{2})^{3}}{3} \xi_{2}
\end{pmatrix}^{3}
\end{aligned}}\\
+ \sqrt[3]{
\begin{aligned}
\begin{pmatrix}
\frac{1+(\tau^{\mathbb{T}}_{2})^{3}+(\tau^{\mathbb{I}}_{2})^{3}-(\tau^{\mathbb{F}}_{2})^{3}}{3} \eta_{2}
- \frac{1+(\tau^{\mathbb{T}}_{3})^{3}+(\tau^{\mathbb{I}}_{3})^{3}-(\tau^{\mathbb{F}}_{3})^{3}}{3} \eta_{3}
\end{pmatrix}^{3}\\
+ \frac{1}{3} 
\begin{pmatrix}
\frac{1+(\tau^{\mathbb{T}}_{2})^{3}+(\tau^{\mathbb{I}}_{2})^{3}-(\tau^{\mathbb{F}}_{2})^{3} }{3} \xi_{2}
- \frac{1+(\tau^{\mathbb{T}}_{3})^{3}+(\tau^{\mathbb{I}}_{3})^{3}- (\tau^{\mathbb{F}}_{3})^{3}}{3} \xi_{3}
\end{pmatrix}^{3}
\end{aligned}}
\end{pmatrix}^{3}\\
&=&\frac{1}{27} 
\begin{bmatrix}
\left((\Phi_{1}\eta_{1}-\Phi_{2}\eta_{2})^{3} + \frac{1}{3}(\Phi_{1}\xi_{1}-\Phi_{2}\xi_{2})^{3}\right)
+\left((\Phi_{2}\eta_{2}-\Phi_{3}\eta_{3})^{3} + \frac{1}{3}(\Phi_{2}\xi_{2}-\Phi_{3}\xi_{3})^{3}\right)\\
+3\left(\sqrt[3]{\big[(\Phi_{1}\eta_{1}-\Phi_{2}\eta_{2})^{3} + \frac{1}{3}(\Phi_{1}\xi_{1}-\Phi_{2}\xi_{2})^{3}\big]^{2}} 
\times \sqrt[3]{ (\Phi_{2}\eta_{2}-\Phi_{3}\eta_{3})^{3} + \frac{1}{3}(\Phi_{2}\xi_{2}-\Phi_{3}\xi_{3})^{3}}\,\,\right)\\
+ 3\left(\sqrt[3]{(\Phi_{1}\eta_{1}-\Phi_{2}\eta_{2})^{3} + \frac{1}{3}(\Phi_{1}\xi_{1}-\Phi_{2}\xi_{2})^{3}} 
 \times 
\sqrt[3]{ \big[(\Phi_{2}\eta_{2}-\Phi_{3}\eta_{3})^{3} + \frac{1}{3}(\Phi_{2}\xi_{2}-\Phi_{3}\xi_{3})^{3}\big]^{2}}\,\,\right)
\end{bmatrix}\\
&=&\frac{1}{27} 
\begin{bmatrix}
\left((\Phi_{1}\eta_{1}-\Phi_{2}\eta_{2})^{3} + \frac{1}{3}(\Phi_{1}\xi_{1}-\Phi_{2}\xi_{2})^{3}\right)
+\left((\Phi_{2}\eta_{2}-\Phi_{3}\eta_{3})^{3} + \frac{1}{3}(\Phi_{2}\xi_{2}-\Phi_{3}\xi_{3})^{3}\right)\\
+3\left(\sqrt[3]{\Big(\big[(\Phi_{1}\eta_{1}-\Phi_{2}\eta_{2})^{3} + \frac{1}{3}(\Phi_{1}\xi_{1}-\Phi_{2}\xi_{2})^{3}\big]^{2}\Big)
\times \Big((\Phi_{2}\eta_{2}-\Phi_{3}\eta_{3})^{3} + \frac{1}{3}(\Phi_{2}\xi_{2}-\Phi_{3}\xi_{3})^{3}\Big)}\,\,\right)\\
+ 3\left(\sqrt[3]{\Big((\Phi_{1}\eta_{1}-\Phi_{2}\eta_{2})^{3} + \frac{1}{3}(\Phi_{1}\xi_{1}-\Phi_{2}\xi_{2})^{3}\Big)
 \times \Big(\big[(\Phi_{2}\eta_{2}-\Phi_{3}\eta_{3})^{3} + \frac{1}{3}(\Phi_{2}\xi_{2}-\Phi_{3}\xi_{3})^{3}\big]^{2}\Big)}\,\,\right)
\end{bmatrix}\\
\end{eqnarray*}
\begin{eqnarray*}
&\geq&\frac{1}{27} 
\begin{bmatrix}
\left((\Phi_{1}\eta_{1}-\Phi_{2}\eta_{2})^{3} + \frac{1}{3}(\Phi_{1}\xi_{1}-\Phi_{2}\xi_{2})^{3}\right)
+\left((\Phi_{2}\eta_{2}-\Phi_{3}\eta_{3})^{3} + \frac{1}{3}(\Phi_{2}\xi_{2}-\Phi_{3}\xi_{3})^{3}\right)\\
+3\left(\sqrt[3]{\Big((\Phi_{1}\eta_{1}-\Phi_{2}\eta_{2})^{2} \cdot (\Phi_{2}\eta_{2}-\Phi_{3}\eta_{3})\Big)^{3} 
 +\frac{1}{27} \Big((\Phi_{1}\xi_{1}-\Phi_{2}\xi_{2})^{2}  \cdot 
 (\Phi_{2}\xi_{2}-\Phi_{3}\xi_{3})\Big)^{3}}\,\right)\\
+3\left(\sqrt[3]{\Big((\Phi_{1}\eta_{1}-\Phi_{2}\eta_{2}) \cdot (\Phi_{2}\eta_{2}-\Phi_{3}\eta_{3})^{2}\Big)^{3} 
 +\frac{1}{27} \Big((\Phi_{1}\xi_{1}-\Phi_{2}\xi_{2}) \cdot 
 (\Phi_{2}\xi_{2}-\Phi_{3}\xi_{3})^{2}\Big)^{3}}\,\right)
\end{bmatrix}\\
&\geq&\frac{1}{27} 
\begin{bmatrix}
\left((\Phi_{1}\eta_{1}-\Phi_{2}\eta_{2})^{3} + \frac{1}{3}(\Phi_{1}\xi_{1}-\Phi_{2}\xi_{2})^{3}\right)
+\left((\Phi_{2}\eta_{2}-\Phi_{3}\eta_{3})^{3} + \frac{1}{3}(\Phi_{2}\xi_{2}-\Phi_{3}\xi_{3})^{3}\right)\\
+3\left(\Big((\Phi_{1}\eta_{1}-\Phi_{2}\eta_{2})^{2} \cdot (\Phi_{2}\eta_{2}-\Phi_{3}\eta_{3})\Big)
 +\frac{1}{3} \Big((\Phi_{1}\xi_{1}-\Phi_{2}\xi_{2})^{2}  \cdot 
 (\Phi_{2}\xi_{2}-\Phi_{3}\xi_{3})\Big)\,\right)\\
+3\left(\Big((\Phi_{1}\eta_{1}-\Phi_{2}\eta_{2}) \cdot (\Phi_{2}\eta_{2}-\Phi_{3}\eta_{3})^{2}\Big)
 +\frac{1}{3} \Big((\Phi_{1}\xi_{1}-\Phi_{2}\xi_{2}) \cdot 
 (\Phi_{2}\xi_{2}-\Phi_{3}\xi_{3})^{2}\Big)\,\right)
\end{bmatrix}\\
&=&\frac{1}{27} 
\begin{bmatrix}
(\Phi_{1}\eta_{1}-\Phi_{2}\eta_{2})^{3} + \frac{1}{3}(\Phi_{1}\xi_{1}-\Phi_{2}\xi_{2})^{3}
+(\Phi_{2}\eta_{2}-\Phi_{3}\eta_{3})^{3} + \frac{1}{3}(\Phi_{2}\xi_{2}-\Phi_{3}\xi_{3})^{3}\\
+3\Big((\Phi_{1}\eta_{1}-\Phi_{2}\eta_{2})^{2} \cdot (\Phi_{2}\eta_{2}-\Phi_{3}\eta_{3})\Big)
+\Big((\Phi_{1}\xi_{1}-\Phi_{2}\xi_{2})^{2}  \cdot 
(\Phi_{2}\xi_{2}-\Phi_{3}\xi_{3})\Big)\\
+3\Big((\Phi_{1}\eta_{1}-\Phi_{2}\eta_{2}) \cdot (\Phi_{2}\eta_{2}-\Phi_{3}\eta_{3})^{2}\Big)
+\Big((\Phi_{1}\xi_{1}-\Phi_{2}\xi_{2}) \cdot 
(\Phi_{2}\xi_{2}-\Phi_{3}\xi_{3})^{2}\Big)
\end{bmatrix}\\
&=& \frac{1}{27}
\begin{bmatrix}
(\Phi_{1}\eta_{1}-\Phi_{2}\eta_{2}+ \Phi_{2}\eta_{2}-\Phi_{3}\eta_{3})^{3}
+\frac{1}{3} (\Phi_{1}\xi_{1}-\Phi_{2}\xi_{2} + \Phi_{2}\xi_{2}-\Phi_{3}\xi_{3})^{3}
\end{bmatrix}\\
&=& \frac{1}{27}
\begin{bmatrix}
(\Phi_{1}\eta_{1}- \Phi_{3}\eta_{3})^{3} +\frac{1}{3} (\Phi_{1}\xi_{1}-\Phi_{3}\xi_{3})^{3}
\end{bmatrix}\\
&=& \mathbb{D}_{E}(L_{1}, L_{3})^{3}
\end{eqnarray*}}
Here, $\Phi_{1}= \frac{1+(\tau^{\mathbb{T}}_{1})^{3}+(\tau^{\mathbb{I}}_{1})^{3}- (\tau^{\mathbb{F}}_{1})^{3}}{3}, \Phi_{2}= \frac{1+(\tau^{\mathbb{T}}_{2})^{3}+ (\tau^{\mathbb{I}}_{2})^{3}-(\tau^{\mathbb{F}}_{2})^{3}}{3}$ and $\Phi_{3}= \frac{1+(\tau^{\mathbb{T}}_{3})^{3}+(\tau^{\mathbb{I}}_{3})^{3}- (\tau^{\mathbb{F}}_{3})^{3}}{3}$.\\
Thus, $ \mathbb{D}_{E}(L_{1}, {L_{3}}) \leq \mathbb{D}_{E}(L_{1}, L_{2}) + \mathbb{D}_{E}(L_{2}, {L_{3}})$.

\begin{thm} 
For any three FNNNs $L_{1} = \big\langle  (\eta_{1}, \xi_{1}); \tau^{\mathbb{T}}_{1},\tau^{\mathbb{I}}_{1}, \tau^{\mathbb{F}}_{1}  \big\rangle$, $L_{2} = \big\langle  (\eta_{2}, \xi_{2}); \tau^{\mathbb{T}}_{2}, \tau^{\mathbb{I}}_{2},\tau^{\mathbb{F}}_{2}  \big\rangle  $ and ${L_{3}} = \big\langle  (\eta_{3}, \xi_{3}); \tau^{\mathbb{T}}_{3}, \tau^{\mathbb{I}}_{3}, \tau^{\mathbb{F}}_{3}  \big\rangle $. Then the following properties are hold:
\begin{enumerate}
\item $\mathbb{D}_{H}(L_{1}, L_{2}) = 0$ if and only if $L_{1}= L_{2}$
\item $\mathbb{D}_{H}(L_{1}, L_{2})= \mathbb{D}_{H}(L_{2}, L_{1})$
\item $ \mathbb{D}_{H}(L_{1}, {L_{3}}) \leq \mathbb{D}_{H}(L_{1}, L_{2}) + \mathbb{D}_{H}(L_{2}, {L_{3}})$.
\end{enumerate}
\end{thm}

\begin{rmk} 
Since $L_{1}=\big\langle \tau^{\mathbb{T}}_{1},\tau^{\mathbb{I}}_{1}, \tau^{\mathbb{F}}_{1} \big\rangle  = \langle 1,1,0 \rangle $, $L_{2} = \big\langle  \tau^{\mathbb{T}}_{2}, \tau^{\mathbb{I}}_{2}, \tau^{\mathbb{F}}_{2} \big\rangle = \langle  1,1,0 \rangle  $, which is known as distance between FNNNs is transformed to the distance between FNNs.
\end{rmk}

\section{Aggregation operators for FNNN} \label{sec:5}
In this section, we will define new aggregation operations for FNNN.

\subsection{FNN weighted averaging (FNNWA) operator}

\begin{defn}
Let $L_{i} = \big\langle  (\eta_{i}, \xi_{i}); \tau^{\mathbb{T}}_{i},\tau^{\mathbb{I}}_{i}, \tau^{\mathbb{F}}_{i}  \big\rangle $ be FNNNs, where $i=1,2,...,n$, $\omega_{i}$ be a weight of $L_{i} $ and $\omega_{i} \geq 0$, $\displaystyle\sum^{n}_{i=1} \omega_{i}=1$. Then, the FNNWA operator can be defined as FNNWA$(L_{1} , L_{2},..., L_{n})= \displaystyle\sum^{n}_{i=1} \omega_{i}L_{i}$.
\end{defn}

\begin{thm}
Let $L_{i} = \big\langle  (\eta_{i}, \xi_{i}); \tau^{\mathbb{T}}_{i},\tau^{\mathbb{I}}_{i} , \tau^{\mathbb{F}}_{i}  \big\rangle  $ be FNNNs, where $i=1,2,...,n$. Then FNNWA operator is defined as\\
FNNWA$(L_{1} , L_{2},..., L_{n})=
\begin{pmatrix}
\Big(\displaystyle\sum^{n}_{i=1} \omega_{i}\eta_{i}, \displaystyle\sum^{n}_{i=1} \omega_{i}\xi_{i}\Big);
 \sqrt[3\Lambda]{1-\displaystyle\prod^{n}_{i=1}\Big(1-(\tau^{\mathbb{T}}_{i})^{3\Lambda}  \Big)^{\omega_{i}}},\\
 \sqrt[\Lambda]{1-\displaystyle\prod^{n}_{i=1}\Big( 1-(\tau^{\mathbb{I}}_{i})^{\Lambda} \Big)^{\omega_{i}}},
 \displaystyle\prod^{n}_{i=1} (\tau^{\mathbb{F}}_{i})^{\omega_{i}}
\end{pmatrix}$.
\end{thm}
\textbf{Proof.} The proof follows from mathematical induction method.\\
If $n=2$, then FNNWA$(L_{1} , L_{2})=  \omega_{1} L_{1} \boxplus \omega_{2}L_{2}$,\\
where,
$$\omega_{1} L_{1}=
\begin{pmatrix}
( \omega_{1}\eta_{1}, \omega_{1}\xi_{1} );
 \sqrt[3\Lambda]{1- \Big( 1-(\tau^{\mathbb{T}}_{1})^{3\Lambda} \Big )^{\omega_{1}}},\\
 \sqrt[\Lambda]{1- \Big( 1-(\tau^{\mathbb{I}}_{1})^{\Lambda} \Big )^{\omega_{1}}},
(\tau^{\mathbb{F}}_{1})^{\omega_{1}}
\end{pmatrix}$$
and
$$\omega_{2} L_{2}=
\begin{pmatrix}
( \omega_{2}\eta_{2}, \omega_{2}\xi_{2} );
 \sqrt[3\Lambda]{1- \Big( 1-(\tau^{\mathbb{T}}_{2})^{3\Lambda}  \Big)^{\omega_{2}}},\\
 \sqrt[\Lambda]{1- \Big( 1-(\tau^{\mathbb{I}}_{2})^{\Lambda} \Big )^{\omega_{2}}},
(\tau^{\mathbb{F}}_{2})^{\omega_{2}}
\end{pmatrix}$$\\
Using Definition \ref{def2}, we get
\begin{eqnarray*}
\omega_{1} L_{1} \boxplus \omega_{2} L_{2}
&=& \begin{pmatrix}
( \omega_{1}\eta_{1} + \omega_{2}\eta_{2}  , \omega_{1}\xi_{1}+ \omega_{2}\xi_{2} );\\
\sqrt[3\Lambda]{
\begin{aligned}
\Big(1- \big( 1-(\tau^{\mathbb{T}}_{1})^{3\Lambda}\big)^{\omega_{1}}\Big) + \Big(1- \big( 1-(\tau^{\mathbb{T}}_{2})^{3\Lambda}\big)^{\omega_{2}}\Big)\\
-\Big(1- \big( 1-(\tau^{\mathbb{T}}_{1})^{3\Lambda}\big)^{\omega_{1}}\Big) \cdot \Big(1- \big( 1-(\tau^{\mathbb{T}}_{2})^{3\Lambda}\big)^{\omega_{2}}\Big),
\end{aligned}}
\,\,\,\, \\
\sqrt[\Lambda]{
\begin{aligned}
\Big(1- \big( 1-(\tau^{\mathbb{I}}_{1})^{\Lambda}\big)^{\omega_{1}}\Big) + \Big(1- \big( 1-(\tau^{\mathbb{I}}_{2})^{\Lambda}\Big)^{\omega_{2}}\Big)\\
-\Big(1- \big( 1-(\tau^{\mathbb{I}}_{1})^{\Lambda}\big)^{\omega_{1}}\Big) \cdot \Big(1- \big( 1-(\tau^{\mathbb{I}}_{2})^{\Lambda}\Big)^{\omega_{2}}\Big),
\end{aligned}}
\,\,\,\, \\
(\tau^{\mathbb{F}}_{1})^{\omega_{1}}\cdot(\tau^{\mathbb{F}}_{2})^{\omega_{2}}
\end{pmatrix}\\
&=&\begin{pmatrix}
( \omega_{1}\eta_{1} + \omega_{2}\eta_{2}  , \omega_{1}\xi_{1}+ \omega_{2}\xi_{2} );\\
\sqrt[3\Lambda]{
\begin{aligned}
1- \Big( 1-(\tau^{\mathbb{T}}_{1})^{3\Lambda}\Big)^{\omega_{1}} \Big( 1-(\tau^{\mathbb{T}}_{2})^{3\Lambda}\Big)^{\omega_{2}},
\end{aligned}}
\,\,\,\\
\sqrt[\Lambda]{
\begin{aligned}
1- \Big( 1-(\tau^{\mathbb{I}}_{1})^{\Lambda}\Big)^{\omega_{1}} \Big( 1-(\tau^{\mathbb{I}}_{2})^{\Lambda}\Big)^{\omega_{2}},
\end{aligned}}
\,\,\,\\
(\tau^{\mathbb{F}}_{1})^{\omega_{1}}\cdot(\tau^{\mathbb{F}}_{2})^{\omega_{2}}
\end{pmatrix}\\
&=&\begin{pmatrix}
\Big(\displaystyle\sum^{2}_{i=1} \omega_{i}\eta_{i}, \displaystyle\sum^{2}_{i=1} \omega_{i}\xi_{i} \Big);
 \sqrt[3\Lambda]{1-\displaystyle\prod^{2}_{i=1}\Big( 1-(\tau^{\mathbb{T}}_{i})^{3\Lambda}  \Big)^{\omega_{i}}},\\
 \sqrt[\Lambda]{1-\displaystyle\prod^{2}_{i=1}\Big( 1-(\tau^{\mathbb{I}}_{i})^{\Lambda}  \Big)^{\omega_{i}}},
 \displaystyle\prod^{2}_{i=1} (\tau^{\mathbb{F}}_{i})^{\omega_{i}}
\end{pmatrix}.
\end{eqnarray*}
Using induction $n \geq 3$, we can show that for $n \geq 3$;\\
FNNWA$(L_{1} , L_{2},..., L_{k})
=\begin{pmatrix}
\Big(\displaystyle\sum^{k}_{i=1} \omega_{i}\eta_{i}, \displaystyle\sum^{k}_{i=1} \omega_{i}\xi_{i} \Big);
 \sqrt[3\Lambda]{1-\displaystyle\prod^{k}_{i=1}\Big( 1-(\tau^{\mathbb{T}}_{i})^{3\Lambda}  \Big)^{\omega_{i}}},\\
 \sqrt[\Lambda]{1-\displaystyle\prod^{k}_{i=1}\Big( 1-(\tau^{\mathbb{I}}_{i})^{\Lambda}  \Big)^{\omega_{i}}},
 \displaystyle\prod^{k}_{i=1} (\tau^{\mathbb{F}}_{i})^{\omega_{i}}
\end{pmatrix}$.\\
If $n= k+1$, then $FNNWA (L_{1} , L_{2},..., L_{k}, L_{k+1})$
\begin{eqnarray*}
&=&\begin{pmatrix}
\left( \displaystyle\sum^{k}_{i=1} \omega_{i}\eta_{i} + \omega_{k+1}\eta_{k+1}  , \displaystyle\sum^{k}_{i=1} \omega_{i}\xi_{i}+ \omega_{k+1}\xi_{k+1} \right);\\
\sqrt[3\Lambda]{
\begin{aligned}
\displaystyle\sum^{k}_{i=1} \Big(1- \Big( 1-(\tau^{\mathbb{T}}_{i})^{3\Lambda}\Big)^{\omega_{i}}\Big) + \Big(1- \Big( 1-(\tau^{\mathbb{T}}_{k+1})^{3\Lambda}\Big)^{\omega_{k+1}}\Big)\\
-\displaystyle\prod^{k}_{i=1} \Big(1- \Big( 1-(\tau^{\mathbb{T}}_{i})^{3\Lambda}\Big)^{\omega_{i}}\Big) \cdot \Big(1- \Big( 1-(\tau^{\mathbb{T}}_{k+1})^{3\Lambda}\Big)^{\omega_{k+1}}\Big),
\end{aligned}}
\,\,\,\, \\
\sqrt[\Lambda]{
\begin{aligned}
\displaystyle\sum^{k}_{i=1} \Big(1- \Big( 1-(\tau^{\mathbb{I}}_{i})^{\Lambda}\Big)^{\omega_{i}}\Big) + \Big(1- \Big( 1-(\tau^{\mathbb{I}}_{k+1})^{\Lambda}\Big)^{\omega_{k+1}}\Big)\\
-\displaystyle\prod^{k}_{i=1}\Big(1- \Big( 1-(\tau^{\mathbb{I}}_{i})^{\Lambda}\Big)^{\omega_{i}}\Big) \cdot \Big(1- \Big( 1-(\tau^{\mathbb{I}}_{k+1})^{\Lambda}\Big)^{\omega_{k+1}}\Big),
\end{aligned}}
\,\,\,\, \\
\displaystyle\prod^{k}_{i=1}(\tau^{\mathbb{F}}_{i})^{\omega_{i}}\cdot(\tau^{\mathbb{F}}_{k+1})^{\omega_{k+1}}
\end{pmatrix}\\
&=& \begin{pmatrix}
\left( \displaystyle\sum^{k+1}_{i=1} \omega_{i}\eta_{i}, \displaystyle\sum^{k+1}_{i=1} \omega_{i}\xi_{i}\right);\\
\sqrt[3\Lambda]{
\begin{aligned}
1-\displaystyle\prod^{k}_{i=1} \Big( 1-(\tau^{\mathbb{T}}_{i})^{3\Lambda}\Big)^{\omega_{i}} \cdot  \Big(1-(\tau^{\mathbb{T}}_{k+1})^{3\Lambda}\Big)^{\omega_{k+1}}\\
\end{aligned}},\\
\sqrt[\Lambda]{
\begin{aligned}
1-\displaystyle\prod^{k}_{i=1} \Big( 1-(\tau^{\mathbb{I}}_{i})^{\Lambda}\Big)^{\omega_{i}} \cdot  \Big(1-(\tau^{\mathbb{I}}_{k+1})^{\Lambda}\Big)^{\omega_{k+1}}
\end{aligned}},
\,\,\,\, \\
\displaystyle\prod^{k+1}_{i=1}(\tau^{\mathbb{F}}_{i})^{\omega_{i}}
\end{pmatrix}\\
&=&\begin{pmatrix}
\left( \displaystyle\sum^{k+1}_{i=1} \omega_{i}\eta_{i}, \displaystyle\sum^{k+1}_{i=1} \omega_{i}\xi_{i}\right);
\sqrt[3\Lambda]{
\begin{aligned}
1-\displaystyle\prod^{k+1}_{i=1} \Big( 1-(\tau^{\mathbb{T}}_{i})^{3\Lambda}\Big)^{\omega_{i}}\\
\end{aligned}},
\,\,\,\, \\
\sqrt[\Lambda]{
\begin{aligned}
1-\displaystyle\prod^{k+1}_{i=1} \Big( 1-(\tau^{\mathbb{I}}_{i})^{\Lambda}\Big)^{\omega_{i}}\\
\end{aligned}},
\,\,\,\,
\displaystyle\prod^{k+1}_{i=1}(\tau^{\mathbb{F}}_{i})^{\omega_{i}}
\end{pmatrix}
\end{eqnarray*}
Hence, this theorem holds for any $k$.

\begin{thm} \label{th5.3}
Let $L_{i} = \big\langle  (\eta_{i}, \xi_{i}); \tau^{\mathbb{T}}_{i}, \tau^{\mathbb{I}}_{i} \tau^{\mathbb{F}}_{i}  \big\rangle$, where $i= 1,2,...,n$ be FNNNs and $L_{i} = L$, then  FNNWA$(L_{1} , L_{2},..., L_{n})= L$.
\end{thm}
\textbf{Proof.}
Given that, $(\eta_{i}, \xi_{i}) = (\eta, \xi)$, $(\tau^{\mathbb{T}}_{i}, \tau^{\mathbb{I}}_{i}, \tau^{\mathbb{F}}_{i})  = (\tau^{\mathbb{T}},\tau^{\mathbb{I}},\tau^{\mathbb{F}}) $, for $i= 1,2,...,n$ and  $\displaystyle\sum^{n}_{i=1}{\omega_{i}}=1$. By using Definition \ref{def2}, we get,
\begin{eqnarray*}
FNNWA (L_{1} , L_{2},..., L_{n})
&=& \begin{pmatrix}
\left( \displaystyle\sum^{n}_{i=1} \omega_{i}\eta_{i}, \displaystyle\sum^{n}_{i=1} \omega_{i}\xi_{i}\right);
\sqrt[3\Lambda]{
\begin{aligned}
1-\displaystyle\prod^{n}_{i=1} \Big( 1-(\tau^{\mathbb{T}})^{3\Lambda}\Big)^{\omega_{i}}\\
\end{aligned}},
\,\,\,\,\\
\sqrt[\Lambda]{
\begin{aligned}
1-\displaystyle\prod^{n}_{i=1} \Big( 1-(\tau^{\mathbb{I}})^{\Lambda}\Big)^{\omega_{i}}\\
\end{aligned}},
\,\,
\displaystyle\prod^{n}_{i=1}(\tau^{\mathbb{F}})^{\omega_{i}}
\end{pmatrix}\\
\end{eqnarray*}
\begin{eqnarray*}
&=&\begin{pmatrix}
\left(\eta \displaystyle\sum^{n}_{i=1} \omega_{i}, \xi\displaystyle\sum^{n}_{i=1} \omega_{i}\right);
\sqrt[3\Lambda]{
\begin{aligned}
1-\Big( 1-(\tau^{\mathbb{T}})^{3\Lambda}\Big)^{\displaystyle\sum^{n}_{i=1} \omega_{i}}\\
\end{aligned}},
\,\, \\
\sqrt[\Lambda]{
\begin{aligned}
1-\Big( 1-(\tau^{\mathbb{I}})^{\Lambda}\Big)^{\displaystyle\sum^{n}_{i=1} \omega_{i}}\\
\end{aligned}},\,\,
(\tau^{\mathbb{F}})^{\displaystyle\sum^{n}_{i=1} \omega_{i}}
\end{pmatrix}\\
&=& \begin{pmatrix}
(\eta, \xi);
\sqrt[3\Lambda]{
\begin{aligned}
1-\Big( 1-(\tau^{\mathbb{T}})^{3\Lambda}\Big)\\
\end{aligned}},
\sqrt[\Lambda]{
\begin{aligned}
1-\Big( 1-(\tau^{\mathbb{I}})^{\Lambda}\Big)\\
\end{aligned}},
(\tau^{\mathbb{F}})
\end{pmatrix}\\
&=& 
\Big((\eta, \xi);
\tau^{\mathbb{T}},\tau^{\mathbb{I}},\tau^{\mathbb{F}}\Big) = L
\end{eqnarray*}

\subsection{FNN weighted geometric(FNNWG) operator}
In this section, we define FNN weighted geometric (FNNWG) operator.

\begin{defn}
Let $L_{i} = \big\langle(\eta_{i}, \xi_{i}); \tau^{\mathbb{T}}_{i},\tau^{\mathbb{I}}_{i}, \tau^{\mathbb{F}}_{i}\big\rangle$ be FNNNs, where $i=1,2,...,n$, $\omega_{i}$ be a weight  of $L_{i}$. Then, the FNNWG operator can be defined as FNNWG $(L_{1}, L_{2},..., L_{n}) = \displaystyle\prod^{n}_{i=1} L_{i}^{\omega_{i}}.$
\end{defn}

\begin{thm}
Let $L_{i} = \big\langle  (\eta_{i}, \xi_{i}); \tau^{\mathbb{T}}_{i}, \tau^{\mathbb{F}}_{i}  \big\rangle  $ be FNNNs where $i=1,2,...,n$. Then, the FNNWG operator can be defined as
$$ FNNWG (L_{1} , L_{2},..., L_{n})=
\begin{pmatrix}
\Big(\displaystyle\prod^{n}_{i=1} \eta_{i}^{\omega_{i}}, \displaystyle\prod^{n}_{i=1} \xi_{i}^{\omega_{i}}\Big );
 \displaystyle\prod^{n}_{i=1} (\tau^{\mathbb{T}}_{i})^{\omega_{i}},\\
 \sqrt[\Lambda]{1-\displaystyle\prod^{n}_{i=1}\Big( 1-(\tau^{\mathbb{I}}_{i})^{\Lambda}  \Big)^{\omega_{i}}},
 \sqrt[3\Lambda]{1-\displaystyle\prod^{n}_{i=1}\Big( 1-(\tau^{\mathbb{F}}_{i})^{3\Lambda} \Big )^{\omega_{i}}}
\end{pmatrix}$$
\end{thm}
\textbf{Proof.} The proof follows from mathematical induction method.\\
If $n=2$, then FNNWG$(L_{1} , L_{2})=   L_{1}^{\omega_{1}} \boxtimes L_{2}^{\omega_{2}}$,\\
where,
$$ L_{1}^{\omega_{1}}=
\begin{pmatrix}
(\eta_{1}^{\omega_{1}}, \xi_{1}^{\omega_{1}} );
(\tau^{\mathbb{T}}_{1})^{\omega_{1}},
 \sqrt[\Lambda]{1- \Big( 1-(\tau^{\mathbb{I}}_{1})^{\Lambda} \Big )^{\omega_{1}}},
 \sqrt[3\Lambda]{1- \Big( 1-(\tau^{\mathbb{F}}_{1})^{3\Lambda} \Big )^{\omega_{1}}}\,\,
\end{pmatrix}$$
$$ L_{2}^{\omega_{2}}=
\begin{pmatrix}
(\eta_{2}^{\omega_{2}}, \xi_{2}^{\omega_{2}} );
(\tau^{\mathbb{T}}_{2})^{\omega_{2}},
 \sqrt[\Lambda]{1- \Big( 1-(\tau^{\mathbb{I}}_{2})^{\Lambda} \Big )^{\omega_{2}}},
 \sqrt[3\Lambda]{1- \Big( 1-(\tau^{\mathbb{F}}_{2})^{3\Lambda}  \Big)^{\omega_{2}}}\,\,
\end{pmatrix}$$\\
By using Definition \ref{def2}, we get,
\begin{eqnarray*}
L_{1}^{\omega_{1}} \boxtimes L_{2}^{\omega_{2}}
&=&
\begin{pmatrix}
\Big( \eta_{1}^{\omega_{1}} \cdot \eta_{2}^{\omega_{2}}, \xi_{1}^{\omega_{1}} \cdot \xi_{2}^{\omega_{2}} \Big);
(\tau^{\mathbb{T}}_{1})^{\omega_{1}} \cdot(\tau^{\mathbb{T}}_{2})^{\omega_{2}},\\
\sqrt[\Lambda]{
\begin{aligned}
\Big(1- \Big( 1-(\tau^{\mathbb{I}}_{1})^{\Lambda}\Big)^{\omega_{1}}\Big) + \Big(1- \Big( 1-(\tau^{\mathbb{I}}_{2})^{\Lambda}\Big)^{\omega_{2}}\Big)\\
-\Big(1- \Big( 1-(\tau^{\mathbb{I}}_{1})^{\Lambda}\Big)^{\omega_{1}}\Big) \cdot \Big(1- \Big( 1-(\tau^{\mathbb{I}}_{2})^{\Lambda}\Big)^{\omega_{2}}\Big),
\end{aligned}}
\,\,\,\, \\
\sqrt[3\Lambda]{
\begin{aligned}
\Big(1- \Big( 1-(\tau^{\mathbb{F}}_{1})^{3\Lambda}\Big)^{\omega_{1}}\Big) + \Big(1- \Big( 1-(\tau^{\mathbb{F}}_{2})^{3\Lambda}\Big)^{\omega_{2}}\Big)\\
-\Big(1- \Big( 1-(\tau^{\mathbb{F}}_{1})^{3\Lambda}\Big)^{\omega_{1}}\Big) \cdot \Big (1-\Big ( 1-(\tau^{\mathbb{F}}_{2})^{3\Lambda}\Big)^{\omega_{2}}\Big)
\end{aligned}}\\
\,\,\,\, 
\end{pmatrix}
\end{eqnarray*}
\begin{eqnarray*}
&=&
\begin{pmatrix}
\Big( \eta_{1}^{\omega_{1}} \cdot \eta_{2}^{\omega_{2}}, \xi_{1}^{\omega_{1}} \cdot \xi_{2}^{\omega_{2}} \Big);
(\tau^{\mathbb{T}}_{1})^{\omega_{1}}\cdot(\tau^{\mathbb{T}}_{2})^{\omega_{2}},\\
\sqrt[\Lambda]{
\begin{aligned}
1- \Big( 1-(\tau^{\mathbb{I}}_{1})^{\Lambda}\Big)^{\omega_{1}}\cdot \Big( 1-(\tau^{\mathbb{I}}_{2})^{\Lambda}\Big)^{\omega_{2}},
\end{aligned}}
\,\,\,\\
\sqrt[3\Lambda]{
\begin{aligned}
1- \Big( 1-(\tau^{\mathbb{F}}_{1})^{3\Lambda}\Big)^{\omega_{1}} \cdot \Big( 1-(\tau^{\mathbb{F}}_{2})^{3\Lambda}\Big)^{\omega_{2}}
\end{aligned}}
\,\,\,
\end{pmatrix}
\end{eqnarray*}\\
Hence, FNNWG$(L_{1} , L_{2})=
\begin{pmatrix}
\Big(\displaystyle\prod^{2}_{i=1} \eta_{i}^{\omega_{i}}, \displaystyle\prod^{2}_{i=1} \xi_{i}^{\omega_{i}} \Big );
 \displaystyle\prod^{2}_{i=1} (\tau^{\mathbb{T}}_{i})^{\omega_{i}},\\
 \sqrt[\Lambda]{1-\displaystyle\prod^{2}_{i=1}\Big( 1-(\tau^{\mathbb{I}}_{i})^{\Lambda} \Big )^{\omega_{i}}},
 \sqrt[3\Lambda]{1-\displaystyle\prod^{2}_{i=1}\Big( 1-(\tau^{\mathbb{F}}_{i})^{3\Lambda} \Big )^{\omega_{i}}}
\end{pmatrix}$\\
Let us consider the $n= k \geq 3$. Then, we have the following.\\
FNNWG$(L_{1} , L_{2},..., L_{k})=
\begin{pmatrix}
\Big(\displaystyle\prod^{k}_{i=1} \eta_{i}^{\omega_{i}}, \displaystyle\prod^{k}_{i=1} \xi_{i}^{\omega_{i}} \Big );
 \displaystyle\prod^{k}_{i=1} (\tau^{\mathbb{T}}_{i})^{\omega_{i}},\\
 \sqrt[\Lambda]{1-\displaystyle\prod^{k}_{i=1}\Big( 1-(\tau^{\mathbb{I}}_{i})^{\Lambda} \Big )^{\omega_{i}}},
 \sqrt[3\Lambda]{1-\displaystyle\prod^{k}_{i=1}\Big( 1-(\tau^{\mathbb{F}}_{i})^{3\Lambda} \Big )^{\omega_{i}}}
\end{pmatrix}$\\
If $n= k+1$, then FNNWG$(L_{1},...,L_{k},L_{k+1})$
\begin{eqnarray*}
&=&\begin{pmatrix}
\Big( \displaystyle\prod^{k}_{i=1} \eta_{i}^{\omega_{i}} \cdot \eta_{k+1}^{^{\omega_{k+1}}}, \displaystyle\prod^{k}_{i=1} \xi_{i}^{^{\omega_{i}}} \cdot \xi_{k+1}^{\omega_{k+1}} \Big);\\
\displaystyle\prod^{k}_{i=1}(\tau^{\mathbb{T}}_{i})^{\omega_{i}}\cdot(\tau^{\mathbb{T}}_{k+1})^{\omega_{k+1}},\\
\sqrt[\Lambda]{
\begin{aligned}
\displaystyle\sum^{k}_{i=1} \Big(1- \Big( 1-(\tau^{\mathbb{I}}_{i})^{\Lambda}\Big)^{\omega_{i}}\Big) + \Big(1-\Big ( 1-(\tau^{\mathbb{I}}_{k+1})^{\Lambda}\Big)^{\omega_{k+1}}\Big)\\
-\displaystyle\prod^{k}_{i=1}\Big(1- \Big( 1-(\tau^{\mathbb{I}}_{i})^{\Lambda}\Big)^{\omega_{i}}\Big) \cdot \Big (1- \Big( 1-(\tau^{\mathbb{I}}_{k+1})^{\Lambda}\Big)^{\omega_{k+1}}\Big)
\end{aligned}},\\
\sqrt[3\Lambda]{
\begin{aligned}
\displaystyle\sum^{k}_{i=1} \Big(1- \Big( 1-(\tau^{\mathbb{F}}_{i})^{3\Lambda}\Big)^{\omega_{i}}\Big) + \Big(1- \Big( 1-(\tau^{\mathbb{F}}_{k+1})^{3\Lambda}\Big)^{\omega_{k+1}}\Big)\\
-\displaystyle\prod^{k}_{i=1}\Big(1- \Big( 1-(\tau^{\mathbb{F}}_{i})^{3\Lambda}\Big)^{\omega_{i}}\Big) \cdot \Big(1- \Big( 1-(\tau^{\mathbb{F}}_{k+1})^{3\Lambda}\Big)^{\omega_{k+1}}\Big)
\end{aligned}}\\
\end{pmatrix}\\
&=& \begin{pmatrix}
\Big( \displaystyle\prod^{k+1}_{i=1} \eta_{i}^{\omega_{i}}, \displaystyle\prod^{k+1}_{i=1} \xi_{i}^{\omega_{i}}\Big);
\displaystyle\prod^{k+1}_{i=1}(\tau^{\mathbb{T}}_{i})^{\omega_{i}},\\
\sqrt[\Lambda]{
\begin{aligned}
1-\displaystyle\prod^{k}_{i=1} \Big( 1-(\tau^{\mathbb{I}}_{i})^{\Lambda}\Big)^{\omega_{i}} \cdot  \Big(1-(\tau^{\mathbb{I}}_{k+1})^{\Lambda}\Big)^{\omega_{k+1}}\\
\end{aligned}},\\
\sqrt[3\Lambda]{
\begin{aligned}
1-\displaystyle\prod^{k}_{i=1} \Big( 1-(\tau^{\mathbb{F}}_{i})^{3\Lambda}\Big)^{\omega_{i}} \cdot  \Big(1-(\tau^{\mathbb{F}}_{k+1})^{3\Lambda}\Big)^{\omega_{k+1}}\\
\end{aligned}}\\
\end{pmatrix}\\
&=& \begin{pmatrix}
\Big( \displaystyle\prod^{k+1}_{i=1} \eta_{i}^{\omega_{i}}, \displaystyle\prod^{k+1}_{i=1} \xi_{i}^{\omega_{i}}\Big);
\displaystyle\prod^{k+1}_{i=1}(\tau^{\mathbb{T}}_{i})^{\omega_{i}},\\
\sqrt[\Lambda]{
\begin{aligned}
1-\displaystyle\prod^{k+1}_{i=1} \Big( 1-(\tau^{\mathbb{I}}_{i})^{\Lambda}\Big)^{\omega_{i}}\\
\end{aligned}},
\sqrt[3\Lambda]{
\begin{aligned}
1-\displaystyle\prod^{k+1}_{i=1} \Big( 1-(\tau^{\mathbb{F}}_{i})^{3\Lambda}\Big)^{\omega_{i}}\\
\end{aligned}}
\end{pmatrix}
\end{eqnarray*}
Therefore, this theorem holds for any natural number.

\begin{thm} 
Let $L_{i} = \big\langle  (\eta_{i}, \xi_{i}); \tau^{\mathbb{T}}_{i}, \tau^{\mathbb{I}}_{i} \tau^{\mathbb{F}}_{i}  \big\rangle$, where $i= 1,2,...,n$ be FNNNs and all are equal with  $L_{i} = L$. Then, FNNWG$(L_{1} , L_{2},..., L_{n})=L$.
\end{thm}
\textbf{Proof.} The proof follows from Theorem \ref{th5.3}.

\subsection{Generalized FNNWA(GFNNWA) operator}
In this section, we define generalized FNNWA operator.

\begin{defn}
Let $L_{i} = \big\langle  (\eta_{i}, \xi_{i}); \tau^{\mathbb{T}}_{i},\tau^{\mathbb{I}}_{i}, \tau^{\mathbb{F}}_{i}  \big\rangle  $ be FNNNs, $\omega_{i}$ be a weight  of $L_{i} $ for $i=1,2,...,n$. Then, GFNNWA$(L_{1} , L_{2},..., L_{n})=
\Big(\displaystyle\sum^{n}_{i=1} \omega_{i}L_{i}^{\Lambda} \Big)^{1/\Lambda}$ is called a generalized FNNWA (GFNNWA) operator.
\end{defn}

\begin{thm}
Let $L_{i} = \big\langle  (\eta_{i}, \xi_{i}); \tau^{\mathbb{T}}_{i},\tau^{\mathbb{I}}_{i} , \tau^{\mathbb{F}}_{i}  \big\rangle  $ be FNNNs, where $i=1,2,...,n$. Then, GFNNWA operator can be defined as \\
GFNNWA $(L_{1} , L_{2},..., L_{n})=
\begin{pmatrix}
\Bigg( \Big(\displaystyle\sum^{n}_{i=1} \omega_{i}\eta_{i}^{\Lambda}\Big)^{1/\Lambda}, \Big(\displaystyle\sum^{n}_{i=1} \omega_{i}\xi_{i}^{\Lambda}\Big)^{1/\Lambda}\Bigg);\\
\Bigg(\sqrt[3\Lambda]{
\begin{aligned}
1-\displaystyle\prod^{n}_{i=1} \Big( 1-\Big((\tau^{\mathbb{T}}_{i})^{\Lambda}\Big)^{3\Lambda}\Big)^{\omega_{i}}\\
\end{aligned}}\,\,\Bigg)^{1/\Lambda},\\
\Bigg(\sqrt[\Lambda]{
\begin{aligned}
1-\displaystyle\prod^{n}_{i=1} \Big( 1-\Big((\tau^{\mathbb{I}}_{i})^{\Lambda}\Big)^{\Lambda}\Big)^{\omega_{i}}\\
\end{aligned}}\,\,\Bigg)^{1/\Lambda},\\
\sqrt[3\Lambda]{1-\Bigg(
\begin{aligned}
1-\Bigg(\displaystyle\prod^{n}_{i=1} \Big( \sqrt[3\Lambda]{1-\big(1-(\tau^{\mathbb{F}}_{i})^{3\Lambda}\big)^{\Lambda}}\,\,\Big)^{\omega_{i}}\\
\end{aligned}\,\,\Bigg)^{3\Lambda}\Bigg)^{1/\Lambda}}\\
\end{pmatrix}$.
\end{thm}
\textbf{Proof.} The proof follows from mathematical induction method. 
First let us show that,
$$\displaystyle\sum^{n}_{i=1} \omega_{i} L_{i}^{\Lambda}=
\begin{pmatrix}
\Bigg( \Big(\displaystyle\sum^{n}_{i=1} \omega_{i}\eta_{i}^{\Lambda}\Big), \Big(\displaystyle\sum^{n}_{i=1} \omega_{i}\xi_{i}^{\Lambda}\Big)\Bigg);
\sqrt[3\Lambda]{
\begin{aligned}
1-\displaystyle\prod^{n}_{i=1} \left( 1-\Big((\tau^{\mathbb{T}}_{i})^{\Lambda}\Big)^{3\Lambda}\right)^{\omega_{i}}\\
\end{aligned}},\\
\sqrt[\Lambda]{
\begin{aligned}
1-\displaystyle\prod^{n}_{i=1} \left( 1-\Big((\tau^{\mathbb{I}}_{i})^{\Lambda}\Big)^{\Lambda}\right)^{\omega_{i}}\\
\end{aligned}},
\begin{aligned}
\displaystyle\prod^{n}_{i=1} \left( \sqrt[3\Lambda]{1-\Big(1-(\tau^{\mathbb{F}}_{i})^{3\Lambda}\Big)^{\Lambda}}\,\,\right)^{\omega_{i}}\\
\end{aligned}
\end{pmatrix}.$$
If $n=2$, then
$\omega_{1} L_{1}^{\Lambda}=
\begin{pmatrix}
\Big(\omega_{1}\eta_{1}^{\Lambda},  \omega_{1}\xi_{1}^{\Lambda}\Big);
\sqrt[3\Lambda]{
\begin{aligned}
1-\bigg( 1-\Big((\tau^{\mathbb{T}}_{1})^{\Lambda}\Big)^{3\Lambda}\bigg)^{\omega_{1}}\\
\end{aligned}} ,\\
\sqrt[\Lambda]{
\begin{aligned}
1-\bigg( 1-\Big((\tau^{\mathbb{I}}_{1})^{\Lambda}\Big)^{\Lambda}\bigg)^{\omega_{1}}\\
\end{aligned}} 
\begin{aligned}
\left( \sqrt[3\Lambda]{1-\Big(1-(\tau^{\mathbb{F}}_{1})^{3\Lambda}\Big)^{\Lambda}}\right)^{\omega_{1}}\\
\end{aligned}
\end{pmatrix}$\\
and \\
$\omega_{2} L_{2}^{\Lambda}=
\begin{pmatrix}
\Big(\omega_{2}\eta_{2}^{\Lambda},  \omega_{2}\xi_{2}^{\Lambda}\Big);
\sqrt[3\Lambda]{
\begin{aligned}
1-\bigg( 1-\Big((\tau^{\mathbb{T}}_{2})^{\Lambda}\Big)^{3\Lambda}\bigg)^{\omega_{1}}\\
\end{aligned}},\\
\sqrt[\Lambda]{
\begin{aligned}
1-\bigg( 1-\Big((\tau^{\mathbb{I}}_{2})^{\Lambda}\Big)^{\Lambda}\bigg)^{\omega_{1}}\\
\end{aligned}},
\begin{aligned}
\left( \sqrt[3\Lambda]{1-\Big(1-(\tau^{\mathbb{F}}_{2})^{3\Lambda}\Big)^{\Lambda}}\right)^{\omega_{1}}\\
\end{aligned}
\end{pmatrix}$.\\\\
By using Definition \ref{def2}, we get, $\omega_{1}L_{1}\boxplus \omega_{2}L_{2}$ can be written as
\begin{eqnarray*}
&=&
\begin{pmatrix}
\Big(\omega_{1}\eta_{1}^{\Lambda}+ \omega_{2}\eta_{2}^{\Lambda}, \omega_{1}\xi_{1}^{\Lambda} + \omega_{2}\xi_{2}^{\Lambda}\Big),\\
\sqrt[3\Lambda]{
\begin{aligned}
\Bigg(\sqrt[3\Lambda]{
\begin{aligned}
1-\bigg( 1-\Big((\tau^{\mathbb{T}}_{1})^{\Lambda}\Big)^{3\Lambda}\bigg)^{\omega_{1}}\\
\end{aligned}}\,\, \Bigg)^{3\Lambda} + \Bigg(\sqrt[3\Lambda]{
\begin{aligned}
1-\bigg( 1-\Big((\tau^{\mathbb{T}}_{2})^{\Lambda}\Big)^{3\Lambda}\bigg)^{\omega_{1}}\\
\end{aligned}}\,\, \Bigg)^{3\Lambda}, \\
- \Bigg(\sqrt[3\Lambda]{
\begin{aligned}
1-\bigg( 1-\Big((\tau^{\mathbb{T}}_{1})^{\Lambda}\Big)^{3\Lambda}\bigg)^{\omega_{1}}\\
\end{aligned}}\,\, \Bigg)^{3\Lambda} \cdot \Bigg(\sqrt[3\Lambda]{
\begin{aligned}
1-\bigg( 1-\Big((\tau^{\mathbb{T}}_{2})^{\Lambda}\Big)^{3\Lambda}\bigg)^{\omega_{1}}\\
\end{aligned}}\,\,\Bigg )^{3\Lambda}
\end{aligned}}\\
\sqrt[\Lambda]{
\begin{aligned}
\Bigg(\sqrt[\Lambda]{
\begin{aligned}
1-\bigg( 1-\Big((\tau^{\mathbb{I}}_{1})^{\Lambda}\Big)^{\Lambda}\bigg)^{\omega_{1}}\\
\end{aligned}}\,\, \Bigg)^{\Lambda} + \Bigg(\sqrt[\Lambda]{
\begin{aligned}
1-\bigg( 1-\Big((\tau^{\mathbb{I}}_{2})^{\Lambda}\Big)^{\Lambda}\bigg)^{\omega_{1}}\\
\end{aligned}}\,\, \Bigg)^{\Lambda}, \\
- \Bigg(\sqrt[\Lambda]{
\begin{aligned}
1-\bigg( 1-\Big((\tau^{\mathbb{I}}_{1})^{\Lambda}\Big)^{\Lambda}\bigg)^{\omega_{1}}\\
\end{aligned}}\,\, \Bigg)^{\Lambda} \cdot \Bigg(\sqrt[\Lambda]{
\begin{aligned}
1-\bigg( 1-\Big((\tau^{\mathbb{I}}_{2})^{\Lambda}\Big)^{\Lambda}\bigg)^{\omega_{1}}\\
\end{aligned}}\,\,\Bigg )^{\Lambda}
\end{aligned}}\\
\begin{aligned}
\left( \sqrt[3\Lambda]{1-\Big(1-(\tau^{\mathbb{F}}_{1})^{3\Lambda}\Big)^{\Lambda}}\right)^{\omega_{1}} \cdot \left( \sqrt[3\Lambda]{1-\Big(1-(\tau^{\mathbb{F}}_{2})^{3\Lambda}\Big)^{\Lambda}}\right)^{\omega_{1}} \\
\end{aligned}\\
\end{pmatrix}
\end{eqnarray*}
\begin{eqnarray*}
&=&
\begin{pmatrix}
\Big(\displaystyle\sum^{2}_{i=1} \omega_{i}\eta_{i}^{\Lambda}, \displaystyle\sum^{2}_{i=1} \omega_{i}\xi_{i}^{\Lambda}\Big),
\sqrt[3\Lambda]{
1- \displaystyle\prod^{2}_{i=1}
\bigg(1-\Big((\tau^{\mathbb{T}}_{1})^{\Lambda}\Big)^{3\Lambda}\bigg)^{\omega_{i}}},\\
\sqrt[\Lambda]{
1- \displaystyle\prod^{2}_{i=1}
\bigg(1-\Big((\tau^{\mathbb{I}}_{1})^{\Lambda}\Big)^{\Lambda}\bigg)^{\omega_{i}}},
\displaystyle\prod^{2}_{i=1}\left( \sqrt[3\Lambda]{1-\Big(1-(\tau^{\mathbb{F}}_{i})^{3\Lambda}\Big)^{\Lambda}}\,\,\right)^{\omega_{i}}
\end{pmatrix}
\end{eqnarray*}
In general,
$$\begin{pmatrix}
\Big(\displaystyle\sum^{k}_{i=1} \omega_{i}\eta_{i}^{\Lambda}, \displaystyle\sum^{k}_{i=1} \omega_{i}\xi_{i}^{\Lambda}\Big);
\sqrt[3\Lambda]{
1- \displaystyle\prod^{k}_{i=1}
\bigg(1-\Big((\tau^{\mathbb{T}}_{1})^{\Lambda}\Big)^{3\Lambda}\bigg)^{\omega_{i}}},\\
\sqrt[\Lambda]{
1- \displaystyle\prod^{k}_{i=1}
\bigg(1-\Big((\tau^{\mathbb{I}}_{1})^{\Lambda}\Big)^{\Lambda}\bigg)^{\omega_{i}}},
\displaystyle\prod^{k}_{i=1}\left( \sqrt[3\Lambda]{1-\Big(1-(\tau^{\mathbb{F}}_{i})^{3\Lambda}\Big)^{\Lambda}}\,\,\right)^{\omega_{i}}
\end{pmatrix}$$.\\
If $n= k+1$, then $\displaystyle\sum^{k}_{i=1} \omega_{i}L_{i}^{\Lambda} +\omega_{k+1} L_{k+1}^{\Lambda}= \displaystyle\sum^{k+1}_{i=1} \omega_{i}L_{i}^{\Lambda}$.\\
Now, $\displaystyle\sum^{k}_{i=1} \omega_{i}L_{i}^{\Lambda} +\omega_{k+1} L_{k+1}^{\Lambda} =  \omega_{1}L_{1}^{\Lambda} \boxplus  \omega_{2}L_{2}^{\Lambda} \boxplus ... \boxplus \omega_{k}L_{k}^{\Lambda} \boxplus \omega_{k+1}L_{k+1}^{\Lambda}$\\
{\small
\begin{eqnarray*}
&=&
\begin{pmatrix}
\Bigg(\displaystyle\sum^{k}_{i=1} \omega_{i}\eta_{i}^{\Lambda}+ \omega_{k+1}\eta_{k+1}^{\Lambda}, \displaystyle\sum^{k}_{i=1} \omega_{i}\xi_{i}^{\Lambda}+ \omega_{k+1}\xi_{k+1}^{\Lambda}\Bigg);\\
\sqrt[3\Lambda]{
\begin{aligned}
\Bigg(\sqrt[3\Lambda]{
\begin{aligned}
1-\displaystyle\prod^{k}_{i=1}\bigg( 1-\Big((\tau^{\mathbb{T}}_{i})^{\Lambda}\Big)^{3\Lambda}\bigg)^{\omega_{i}}\\
\end{aligned}}\,\,\Bigg )^{3\Lambda} + \Bigg(\sqrt[3\Lambda]{
\begin{aligned}
1-\bigg( 1-\Big((\tau^{\mathbb{T}}_{k+1})^{\Lambda}\Big)^{3\Lambda}\bigg)^{\omega_{1}}\\
\end{aligned}}\,\, \Bigg)^{3\Lambda}, \\
- \Bigg(\sqrt[3\Lambda]{
\begin{aligned}
1-\displaystyle\prod^{k}_{i=1}\bigg( 1-\Big((\tau^{\mathbb{T}}_{i})^{\Lambda}\Big)^{3\Lambda}\bigg)^{\omega_{i}}\\
\end{aligned}}\,\,\Bigg )^{3\Lambda} \cdot \Bigg(\sqrt[3\Lambda]{
\begin{aligned}
1-\bigg( 1-\Big((\tau^{\mathbb{T}}_{k+1})^{\Lambda}\Big)^{3\Lambda}\bigg)^{\omega_{1}}\\
\end{aligned}}\,\, \Bigg)^{3\Lambda}
\end{aligned}}\\
\sqrt[\Lambda]{
\begin{aligned}
\Bigg(\sqrt[\Lambda]{
\begin{aligned}
1-\displaystyle\prod^{k}_{i=1}\bigg( 1-\Big((\tau^{\mathbb{I}}_{i})^{\Lambda}\Big)^{\Lambda}\bigg)^{\omega_{i}}\\
\end{aligned}}\,\,\Bigg )^{\Lambda} + \Bigg(\sqrt[\Lambda]{
\begin{aligned}
1-\bigg( 1-\Big((\tau^{\mathbb{I}}_{k+1})^{\Lambda}\Big)^{\Lambda}\bigg)^{\omega_{1}}\\
\end{aligned}}\,\, \Bigg)^{\Lambda}, \\
- \Bigg(\sqrt[\Lambda]{
\begin{aligned}
1-\displaystyle\prod^{k}_{i=1}\bigg( 1-\Big((\tau^{\mathbb{I}}_{i})^{\Lambda}\Big)^{\Lambda}\bigg)^{\omega_{i}}\\
\end{aligned}}\,\,\Bigg )^{\Lambda} \cdot \Bigg(\sqrt[\Lambda]{
\begin{aligned}
1-\bigg( 1-\Big((\tau^{\mathbb{I}}_{k+1})^{\Lambda}\Big)^{\Lambda}\bigg)^{\omega_{1}}\\
\end{aligned}}\,\, \Bigg)^{\Lambda}
\end{aligned}}\\
\begin{aligned}
\displaystyle\prod^{k}_{i=1}\left( \sqrt[3\Lambda]{1-\Big(1-(\tau^{\mathbb{F}}_{i})^{3\Lambda}\Big)^{\Lambda}}\right)^{\omega_{i}} \cdot \left( \sqrt[3\Lambda]{1-\Big(1-(\tau^{\mathbb{F}}_{k+1})^{3\Lambda}\Big)^{\Lambda}}\right)^{\omega_{1}} \\
\end{aligned}\\
\end{pmatrix}
\end{eqnarray*}}
\begin{eqnarray*}
&=&\begin{pmatrix}
\Bigg(\displaystyle\sum^{k+1}_{i=1} \omega_{i}\eta_{i}^{\Lambda}, \displaystyle\sum^{k+1}_{i=1} \omega_{i}\xi_{i}^{\Lambda}\Bigg);
\sqrt[3\Lambda]{
1- \displaystyle\prod^{k+1}_{i=1}
\bigg(1-\Big((\tau^{\mathbb{T}}_{1})^{\Lambda}\Big)^{3\Lambda}\bigg)^{\omega_{i}}},\\
\sqrt[\Lambda]{
1- \displaystyle\prod^{k+1}_{i=1}
\bigg(1-\Big((\tau^{\mathbb{I}}_{1})^{\Lambda}\Big)^{\Lambda}\bigg)^{\omega_{i}}},
\displaystyle\prod^{k+1}_{i=1}\left( \sqrt[3\Lambda]{1-\Big(1-(\tau^{\mathbb{F}}_{i})^{3\Lambda}\Big)^{\Lambda}}\right)^{\omega_{i}}
\end{pmatrix}\\
\end{eqnarray*}
and
$\displaystyle\sum^{k+1}_{i=1} \Big(\omega_{i}L_{i}^{\Lambda}\Big)^{^{1/\Lambda}}=\begin{pmatrix}
\Bigg( \Big(\displaystyle\sum^{k+1}_{i=1} \omega_{i}\eta_{i}^{\Lambda}\Big)^{1/\Lambda}, \Big(\displaystyle\sum^{k+1}_{i=1} \omega_{i}\xi_{i}^{\Lambda}\Big)^{1/\Lambda}\Bigg);\\
\Bigg(\sqrt[3\Lambda]{
\begin{aligned}
1-\displaystyle\prod^{k+1}_{i=1} \Big( 1-\Big((\tau^{\mathbb{T}}_{i})^{\Lambda}\Big)^{3\Lambda}\Big)^{\omega_{i}}\\
\end{aligned}}\,\,\Bigg)^{1/\Lambda},\\
\Bigg(\sqrt[\Lambda]{
\begin{aligned}
1-\displaystyle\prod^{k+1}_{i=1} \Big( 1-\Big((\tau^{\mathbb{I}}_{i})^{\Lambda}\Big)^{\Lambda}\Big)^{\omega_{i}}\\
\end{aligned}}\,\,\Bigg)^{1/\Lambda},\\
\sqrt[3\Lambda]{1-\Bigg(
\begin{aligned}
1-\Bigg(\displaystyle\prod^{k+1}_{i=1} \Bigg( \sqrt[3\Lambda]{1-\Big(1-(\tau^{\mathbb{F}}_{i})^{3\Lambda}\Big)^{\Lambda}}\Bigg)^{\omega_{i}}\\
\end{aligned}\,\,\Bigg)^{3}\Bigg)^{1/\Lambda}}\\
\end{pmatrix}$.\\
This theorem holds for any $k\in \mathbb{N}$. Thus, theorem is proved.

\begin{rmk}
If $\Lambda=1$, then GFNNWA operator is reduced to the FNNWA operator. 
\end{rmk}

\begin{thm} 
Let $L_{i} = \big\langle  (\eta_{i}, \xi_{i}); \tau^{\mathbb{T}}_{i}, \tau^{\mathbb{I}}_{i} \tau^{\mathbb{F}}_{i}  \big\rangle $, where $i= 1,2,...,n$ be FNNNs and all are equal with $L$, that is $L_{i} = L$. Then, GFNNWA $(L_{1} , L_{2},..., L_{n})=L$.
\end{thm}
\textbf{Proof.} The proof follows from Theorem \ref{th5.3}.

\subsection{Generalized FNNWG (GFNNWG) operator}
In this section, we define generalized FNNWG operator.

\begin{defn}
Let $L_{i} = \big\langle  (\eta_{i}, \xi_{i}); \tau^{\mathbb{T}}_{i}, \tau^{\mathbb{I}}_{i}, \tau^{\mathbb{F}}_{i}  \big\rangle  $ be FNNNs, $\omega_{i})$ be a weight  of $L_{i} $, where $i=1,2,...,n$. Then, GFNNWG$(L_{1} , L_{2},..., L_{n})= \frac{1}{\Lambda}\Big(\displaystyle\prod^{n}_{i=1} (\Lambda L_{i})^{\omega_{i}} \Big)$ is called a GFNNWG operator.
\end{defn}

\begin{thm}
Let $L_{i} = \big\langle  (\eta_{i}, \xi_{i}); \tau^{\mathbb{T}}_{i}, \tau^{\mathbb{I}}_{i}, \tau^{\mathbb{F}}_{i}  \big\rangle  $ be the collection of FNNNs. Then GFNNWG operator can be defined as\\
GFNNWG $(L_{1} , L_{2},..., L_{n})=
\begin{pmatrix}
\Bigg( \frac{1}{\Lambda} \displaystyle\prod^{n}_{i=1} (\Lambda \eta_{i})^{\omega_{i}}, \frac{1}{\Lambda} \displaystyle\prod^{n}_{i=1} (\Lambda \xi_{i})^{\omega_{i}}\Bigg);\\
\sqrt[3\Lambda]{1-\Bigg(
\begin{aligned}
1-\Bigg(\displaystyle\prod^{n}_{i=1} \left( \sqrt[3\Lambda]{1-\Big(1-(\tau^{\mathbb{T}}_{i})^{3\Lambda}\Big)^{\Lambda}}\right)^{\omega_{i}}\\
\end{aligned}\,\,\Bigg)^{3\Lambda}\Bigg)^{1/\Lambda}},\\
\Bigg(\sqrt[\Lambda]{
\begin{aligned}
1-\displaystyle\prod^{n}_{i=1} \left( 1-\Big((\tau^{\mathbb{I}}_{i})^{\Lambda}\Big)^{\Lambda}\right)^{\omega_{i}}\\
\end{aligned}}\,\,\Bigg)^{1/\Lambda},\\
\Bigg(\sqrt[3\Lambda]{
\begin{aligned}
1-\displaystyle\prod^{n}_{i=1} \left( 1-\Big((\tau^{\mathbb{F}}_{i})^{\Lambda}\Big)^{3\Lambda}\right)^{\omega_{i}}\\
\end{aligned}}\,\,\Bigg)^{1/\Lambda}\,\, 
\end{pmatrix}$
\end{thm}
\textbf{Proof.} The proof follows from mathematical induction. First let us show that
$$\displaystyle\prod^{n}_{i=1}(\Lambda L_{i})^{\omega_{i}}=
\begin{pmatrix}
\Bigg( \displaystyle\prod^{n}_{i=1} (\Lambda\eta_{i})^{\omega_{i}}, \displaystyle\prod^{n}_{i=1} (\Lambda\xi_{i})^{\omega_{i}}\Bigg);
\begin{aligned}
\displaystyle\prod^{n}_{i=1} \left( \sqrt[3\Lambda]{1-\Big(1-(\tau^{\mathbb{T}}_{i})^{3\Lambda}\Big)^{\Lambda}}\right)^{\omega_{i}}\\
\end{aligned},\\
\sqrt[\Lambda]{
\begin{aligned}
1-\displaystyle\prod^{n}_{i=1} \Big( 1-\Big((\tau^{\mathbb{I}}_{i})^{\Lambda}\Big)^{\Lambda}\Big)^{\omega_{i}}\\
\end{aligned}}, 
\sqrt[3\Lambda]{
\begin{aligned}
1-\displaystyle\prod^{n}_{i=1} \Big( 1-\Big((\tau^{\mathbb{F}}_{i})^{\Lambda}\Big)^{3\Lambda}\Big)^{\omega_{i}}\\
\end{aligned}} 
\end{pmatrix}$$
If $n=2$, then
$$ (\Lambda L_{1})^{\omega_{1}}=
\begin{pmatrix}
\Big((\Lambda\eta_{1})^{\omega_{1}},  (\Lambda\xi_{1})^{\omega_{1}}\Big);
\begin{aligned}
\left( \sqrt[3\Lambda]{1-\Big(1-(\tau^{\mathbb{T}}_{1})^{3\Lambda}\Big)^{\Lambda}}\right)^{\omega_{1}}\\
\end{aligned},\\
\sqrt[\Lambda]{
\begin{aligned}
1-\bigg( 1-\Big((\tau^{\mathbb{I}}_{1})^{\Lambda}\Big)^{\Lambda}\bigg)^{\omega_{1}}\\
\end{aligned}}\,\,
\sqrt[3\Lambda]{
\begin{aligned}
1-\bigg( 1-\Big((\tau^{\mathbb{F}}_{1})^{\Lambda}\Big)^{3\Lambda}\bigg)^{\omega_{1}}\\
\end{aligned}}\,\,
\end{pmatrix}$$ and 
$$ (\Lambda L_{2})^{\omega_{2}}=
\begin{pmatrix}
\Big((\Lambda\eta_{2})^{\omega_{2}}, (\Lambda\xi_{2})^{\omega_{2}}\Big);
\begin{aligned}
\left( \sqrt[3\Lambda]{1-\Big(1-(\tau^{\mathbb{T}}_{2})^{3\Lambda}\Big)^{\Lambda}}\right)^{\omega_{1}}\\
\end{aligned},\\
\sqrt[\Lambda]{
\begin{aligned}
1-\bigg( 1-\Big((\tau^{\mathbb{I}}_{2})^{\Lambda}\Big)^{\Lambda}\bigg)^{\omega_{1}}\\
\end{aligned}}
\sqrt[3\Lambda]{
\begin{aligned}
1-\bigg( 1-\Big((\tau^{\mathbb{F}}_{2})^{\Lambda}\Big)^{3\Lambda}\bigg)^{\omega_{1}}\\
\end{aligned}}
\end{pmatrix}$$
By using Definition \ref{def2}, we get,
$(\Lambda L_{1})^{\omega_{1}}\boxtimes (\Lambda L_{2})^{\omega_{2}}$
\begin{eqnarray*}
&=&\begin{pmatrix}
\Big((\Lambda\eta_{1})^{\omega_{1}} \cdot(\Lambda\eta_{2})^{\omega_{2}}, (\Lambda\xi_{1})^{\omega_{1}} \cdot (\Lambda\xi_{2})^{\omega_{2}}\Big);\\
\begin{aligned}
\Bigg( \sqrt[3\Lambda]{1-\Big(1-(\tau^{\mathbb{T}}_{1})^{3\Lambda}\Big)^{\Lambda}}\Bigg)^{\omega_{1}} \cdot \Bigg( \sqrt[3\Lambda]{1-\Big(1-(\tau^{\mathbb{T}}_{2})^{3\Lambda}\Big)^{\Lambda}}\Bigg)^{\omega_{1}}
\end{aligned},\\
\sqrt[\Lambda]{
\begin{aligned}
\Bigg(\sqrt[\Lambda]{
\begin{aligned}
1-\bigg( 1-\Big((\tau^{\mathbb{I}}_{1})^{\Lambda}\Big)^{\Lambda}\bigg)^{\omega_{1}}\\
\end{aligned}}\,\,\Bigg )^{\Lambda} + \Bigg(\sqrt[\Lambda]{
\begin{aligned}
1-\bigg( 1-\Big((\tau^{\mathbb{I}}_{2})^{\Lambda}\Big)^{\Lambda}\bigg)^{\omega_{1}}\\
\end{aligned}}\,\, \Bigg)^{\Lambda} \\
- \Bigg(\sqrt[\Lambda]{
\begin{aligned}
1-\bigg( 1-\Big((\tau^{\mathbb{I}}_{1})^{\Lambda}\Big)^{\Lambda}\bigg)^{\omega_{1}}\\
\end{aligned}}\,\,\Bigg )^{\Lambda} \cdot \Bigg(\sqrt[\Lambda]{
\begin{aligned}
1-\bigg( 1-\Big((\tau^{\mathbb{I}}_{2})^{\Lambda}\Big)^{\Lambda}\bigg)^{\omega_{1}}\\
\end{aligned}}\,\,\Bigg )^{\Lambda}
\end{aligned}},\\
\sqrt[3\Lambda]{
\begin{aligned}
\Bigg(\sqrt[3\Lambda]{
\begin{aligned}
1-\bigg( 1-\Big((\tau^{\mathbb{F}}_{1})^{\Lambda}\Big)^{3\Lambda}\bigg)^{\omega_{1}}\\
\end{aligned}}\,\,\Bigg )^{3\Lambda} + \Bigg(\sqrt[3\Lambda]{
\begin{aligned}
1-\bigg( 1-\Big((\tau^{\mathbb{F}}_{2})^{\Lambda}\Big)^{3\Lambda}\bigg)^{\omega_{1}}\\
\end{aligned}}\,\, \Bigg)^{3\Lambda} \\
- \Bigg(\sqrt[3\Lambda]{
\begin{aligned}
1-\bigg( 1-\Big((\tau^{\mathbb{F}}_{1})^{\Lambda}\Big)^{3\Lambda}\bigg)^{\omega_{1}}\\
\end{aligned}}\,\,\Bigg )^{3\Lambda} \cdot \Bigg(\sqrt[3\Lambda]{
\begin{aligned}
1-\bigg( 1-\Big((\tau^{\mathbb{F}}_{2})^{\Lambda}\Big)^{3\Lambda}\bigg)^{\omega_{1}}\\
\end{aligned}}\,\,\Bigg )^{3\Lambda}
\end{aligned}}
\end{pmatrix}\\
&=&\begin{pmatrix}
\Bigg(\displaystyle\prod^{2}_{i=1} (\Lambda\eta_{i})^{\omega_{i}}, \displaystyle\prod^{2}_{i=1} (\Lambda\xi_{i})^{\omega_{i}}\Bigg);
\displaystyle\prod^{2}_{i=1}\left(\sqrt[3\Lambda]{1-\Big(1-(\tau^{\mathbb{T}}_{i})^{3\Lambda}\Big)^{\Lambda}}\right)^{\omega_{i}},\\
\sqrt[\Lambda]{
1- \displaystyle\prod^{2}_{i=1}
\bigg(1-\Big((\tau^{\mathbb{I}}_{i})^{\Lambda}\Big)^{\Lambda}\bigg)^{\omega_{i}}},
\sqrt[3\Lambda]{
1- \displaystyle\prod^{2}_{i=1}
\bigg(1-\Big((\tau^{\mathbb{F}}_{i})^{\Lambda}\Big)^{3\Lambda}\bigg)^{\omega_{i}}}
\end{pmatrix}
\end{eqnarray*}
If $n=k$, then
$$\displaystyle\prod^{k}_{i=1} (\Lambda\eta_{i})^{ \omega_{i}}
=\begin{pmatrix}
\Bigg(\displaystyle\prod^{k}_{i=1} (\Lambda\eta_{i})^{\omega_{i}}, \displaystyle\prod^{k}_{i=1} (\Lambda\xi_{i})^{\omega_{i}}\Bigg);
\displaystyle\prod^{k}_{i=1}\left(\sqrt[3\Lambda]{1-\Big(1-(\tau^{\mathbb{T}}_{i})^{3\Lambda}\Big)^{\Lambda}}\right)^{\omega_{i}},\\
\sqrt[\Lambda]{
1- \displaystyle\prod^{k}_{i=1}
\bigg(1-\Big((\tau^{\mathbb{I}}_{i})^{\Lambda}\Big)^{\Lambda}\bigg)^{\omega_{i}}},
\sqrt[3\Lambda]{
1- \displaystyle\prod^{k}_{i=1}
\bigg(1-\Big((\tau^{\mathbb{F}}_{i})^{\Lambda}\Big)^{3\Lambda}\bigg)^{\omega_{i}}}
\end{pmatrix}$$
If $n= k+1$, then $\displaystyle\prod^{k}_{i=1} (\Lambda L_{i})^{\omega_{i}} \cdot (\Lambda L_{k+1})^{\omega_{k+1}}= \displaystyle\prod^{k+1}_{i=1} (\Lambda L_{i})^{\omega_{i}}$.\\
Now, $\displaystyle\prod^{k}_{i=1} (\Lambda L_{i})^{\omega_{i}} \cdot (\Lambda L_{k+1})^{\omega_{k+1}} =  (\Lambda L_{1})^{\omega_{1}} \boxtimes (\Lambda L_{2})^{\omega_{2}} \boxtimes ... \boxtimes (\Lambda L_{k})^{\omega_{k}} \boxtimes (\Lambda L_{k+1})^{\omega_{k+1}}$
\begin{eqnarray*}
&=&\begin{pmatrix}
\Big(\displaystyle\prod^{k}_{i=1} (\Lambda\eta_{i})^{\omega_{i}} \cdot (\Lambda\eta_{k+1})^{\omega_{k+1}}, \displaystyle\prod^{k}_{i=1} (\Lambda\xi_{i})^{\omega_{i}} \cdot (\Lambda\xi_{k+1})^{\omega_{k+1}}\Big);\\
\begin{aligned}
\displaystyle\prod^{k}_{i=1}\Big( \sqrt[3\Lambda]{1-\Big(1-(\tau^{\mathbb{T}}_{i})^{3\Lambda}\Big)^{\Lambda}}\Big)^{\omega_{i}} \cdot \Big( \sqrt[3\Lambda]{1-\Big(1-(\tau^{\mathbb{T}}_{k+1})^{3\Lambda}\Big)^{\Lambda}}\Big)^{\omega_{1}} \\
\end{aligned},\\
\sqrt[\Lambda]{
\begin{aligned}
\Bigg(\sqrt[\Lambda]{
\begin{aligned}
1-\displaystyle\prod^{k}_{i=1}\bigg( 1-\Big((\tau^{\mathbb{I}}_{i})^{\Lambda}\Big)^{\Lambda}\bigg)^{\omega_{i}}\\
\end{aligned}}\,\,\Bigg )^{\Lambda} + \Bigg(\sqrt[\Lambda]{
\begin{aligned}
1-\bigg( 1-\Big((\tau^{\mathbb{I}}_{k+1})^{\Lambda}\Big)^{\Lambda}\bigg)^{\omega_{1}}\\
\end{aligned}}\,\, \Bigg)^{\Lambda}\\
- \Bigg(\sqrt[\Lambda]{
\begin{aligned}
1-\displaystyle\prod^{k}_{i=1}\bigg( 1-\Big((\tau^{\mathbb{I}}_{i})^{\Lambda}\Big)^{\Lambda}\bigg)^{\omega_{i}}\\
\end{aligned}}\,\, \Bigg)^{\Lambda} \cdot \Bigg(\sqrt[\Lambda]{
\begin{aligned}
1-\bigg( 1-\Big((\tau^{\mathbb{I}}_{k+1})^{\Lambda}\Big)^{\Lambda}\bigg)^{\omega_{1}}\\
\end{aligned}}\,\,\Bigg )^{\Lambda}
\end{aligned}},\\
\sqrt[3\Lambda]{
\begin{aligned}
\Bigg(\sqrt[3\Lambda]{
\begin{aligned}
1-\displaystyle\prod^{k}_{i=1}\bigg( 1-\Big((\tau^{\mathbb{F}}_{i})^{\Lambda}\Big)^{3\Lambda}\bigg)^{\omega_{i}}\\
\end{aligned}}\,\,\Bigg )^{3\Lambda} + \Bigg(\sqrt[3\Lambda]{
\begin{aligned}
1-\bigg( 1-\Big((\tau^{\mathbb{F}}_{k+1})^{\Lambda}\Big)^{3\Lambda}\bigg)^{\omega_{1}}\\
\end{aligned}}\,\, \Bigg)^{3\Lambda}\\
- \Bigg(\sqrt[3\Lambda]{
\begin{aligned}
1-\displaystyle\prod^{k}_{i=1}\bigg( 1-\Big((\tau^{\mathbb{F}}_{i})^{\Lambda}\Big)^{3\Lambda}\bigg)^{\omega_{i}}\\
\end{aligned}}\,\, \Bigg)^{3\Lambda} \cdot \Bigg(\sqrt[3\Lambda]{
\begin{aligned}
1-\bigg( 1-\Big((\tau^{\mathbb{F}}_{k+1})^{\Lambda}\Big)^{3\Lambda}\bigg)^{\omega_{1}}\\
\end{aligned}}\,\,\Bigg )^{3\Lambda}
\end{aligned}}\\
\end{pmatrix}\\
&=&\begin{pmatrix}
\Bigg(\displaystyle\prod^{k+1}_{i=1} (\Lambda\eta_{i})^{\omega_{i}}, \displaystyle\prod^{k+1}_{i=1} (\Lambda\xi_{i})^{\omega_{i}}\Bigg),
\displaystyle\prod^{k+1}_{i=1}\Bigg(\sqrt[3\Lambda]{1-\Big(1-(\tau^{\mathbb{T}}_{i})^{3\Lambda}\Big)^{\Lambda}}\Bigg)^{\omega_{i}},\\
\sqrt[\Lambda]{1- \displaystyle\prod^{k+1}_{i=1}
\bigg(1-\Big((\tau^{\mathbb{I}}_{1})^{\Lambda}\Big)^{\Lambda}\bigg)^{\omega_{i}}},
\sqrt[3\Lambda]{1- \displaystyle\prod^{k+1}_{i=1}
\bigg(1-\Big((\tau^{\mathbb{F}}_{1})^{\Lambda}\Big)^{3\Lambda}\bigg)^{\omega_{i}}}
\end{pmatrix}
\end{eqnarray*}

\begin{eqnarray*}
\frac{1}{\Lambda}\Big(\displaystyle\prod^{k+1}_{i=1} (\Lambda L_{i})^{\omega_{i}}\Big)
&=& \begin{pmatrix}
\Bigg( \frac{1}{\Lambda}\Big(\displaystyle\prod^{k+1}_{i=1} (\Lambda\eta_{i})^{\omega_{i}}\Big), \frac{1}{\Lambda}\Big(\displaystyle\prod^{k+1}_{i=1} (\Lambda\xi_{i})^{\omega_{i}}\Big)\Bigg);\\
\sqrt[3\Lambda]{1-\Bigg(
\begin{aligned}
1-\Bigg(\displaystyle\prod^{k+1}_{i=1} \left( \sqrt[3\Lambda]{1-\Big(1-(\tau^{\mathbb{T}}_{i})^{3\Lambda}\Big)^{\Lambda}}\right)^{\omega_{i}}\\
\end{aligned}\,\,\Bigg)^{3\Lambda}\Bigg)^{1/\Lambda}},\\
\Bigg(\sqrt[\Lambda]{
\begin{aligned}
1-\displaystyle\prod^{k+1}_{i=1} \Big( 1-\Big((\tau^{\mathbb{I}}_{i})^{\Lambda}\Big)^{\Lambda}\Big)^{\omega_{i}}\\
\end{aligned}}\,\,\Bigg)^{1/\Lambda},\\
\Bigg(\sqrt[3\Lambda]{
\begin{aligned}
1-\displaystyle\prod^{k+1}_{i=1} \Big( 1-\Big((\tau^{\mathbb{F}}_{i})^{\Lambda}\Big)^{3\Lambda}\Big)^{\omega_{i}}\\
\end{aligned}}\,\,\Bigg)^{1/\Lambda}
\end{pmatrix}
\end{eqnarray*}
This theorem holds for for any $k\in \mathbb{N}$. Thus, the theorem is proved.

\begin{rmk}
If $\Lambda=1$, then GFNNWG operator is reduced to the FNNWG operator. 
\end{rmk}

\begin{thm} 
Let $L_{i} = \big\langle  (\eta_{i}, \xi_{i}); \tau^{\mathbb{T}}_{i}, \tau^{\mathbb{I}}_{i} \tau^{\mathbb{F}}_{i}  \big\rangle  \,\,, (i= 1,2,...,n)$ be the collection of FNNNs and all are equal with  $L_{i} = L$, then GFNNWG $(L_{1} , L_{2},..., L_{n})= L$.
\end{thm}
\textbf{Proof.} The proof follows from Theorem \ref{th5.3}.

\section{MADM based on FNN information} \label{sec:6}
In this section, of multi attribute decision making based on FNN information.\\
Let $\mathbb{E}= \{\mathbb{E}_{1},\mathbb{E}_{2},...,\mathbb{E}_{n}\}$ represents the set of $n$ alternatives, $\mathbb{A}= \{e_{1},e_{2},...,e_{m}\}$ represent set of $m$ attributes, $w = \{\omega_{1},\omega_{2},...,\omega_{m}\}$ represents weights of the set of attributes, for  $i= 1,2,...,n$ and $j= 1,2,...,m$ $\mathbb{E}_{ij} = \big\langle  (\eta_{ij}, \xi_{ij}); \tau^{\mathbb{T}}_{ij}, \tau^{\mathbb{I}}_{ij} \tau^{\mathbb{F}}_{ij} \big\rangle $ is a FNNN of alternative $\mathbb{E}_{i}$ in attribute $e_{j}$. We consider $\tau^{\mathbb{T}}_{ij},\tau^{\mathbb{I}}_{ij}, \tau_{ij}^{\mathbb{F}}\in [0,1]$ and $0 \leq (\tau^{\mathbb{T}}_{ij}(u))^{3}+(\tau^{\mathbb{I}}_{ij}(u))^{3} +(\tau_{ij}^{\mathbb{F}}(u))^{3} \leq 2$. Here $n$-alternative sets and $m$-attribute sets gives $n \times m$ decision matrix. We denote the decision matrix by $\mathbb{D}=(\mathbb{E}_{ij})_{n \times m}$.
In an MADM problem, one tries to make the best possible decision out of finite alternatives, by taking into a collection of attributes with preferred weights. Here all the alternatives against each attributes are represented by FNNNs and the concept of euclidean and hamming distance measures are used to obtain a decision. Summarizing the positive and negative ideal values of each attributes concerning every alternatives the representation is obtained. A decision is achieved by implementing the following algorithm.

\subsection{Algorithm}
In this section, we discuss the algorithm for FNN information. The algorithm for the selection of the optimal choice is given below.\\
{\bf Step1:} Input the decision values for each alternatives.\\
{\bf Step2:} Determine the normalize decision values for each alternatives. The decision matrix $\mathbb{D}=(\mathbb{E}_{ij})_{n \times m}$ is normalized into $\overrightarrow{\mathbb{D}}=(\overbracket{\mathbb{E}_{ij}})_{n \times m}$; where 
$\overbracket{\mathbb{E}_{ij}}= \Big\langle  (\overrightarrow{\eta_{ij}}, \overrightarrow{\xi_{ij}}); \overrightarrow{\tau^{\mathbb{T}}_{ij}}, \overrightarrow{\tau^{\mathbb{I}}_{ij}},\overrightarrow{\tau^{\mathbb{F}}_{ij}} \Big\rangle $ and
$$\overrightarrow{\eta_{ij}}= \frac{\eta_{ij}} {\max_{i}(\eta_{ij})}, \overrightarrow{\xi_{ij}}= \frac{\xi_{ij}}{\max_{i}(\xi_{ij})} \cdot \frac{\xi_{ij}}{\eta_{ij}},\,\, \overrightarrow{\tau^{\mathbb{T}}_{ij}}= \tau^{\mathbb{T}}_{ij}, \overrightarrow{\tau^{\mathbb{I}}_{ij}}= \tau^{\mathbb{I}}_{ij}, \overrightarrow{\tau^{\mathbb{F}}_{ij}}= \tau^{\mathbb{F}}_{ij}.$$
{\bf Step3:} Find the aggregate of the values of each alternative. On the basis of FNN information aggregation operators and attribute $e_{j}$ in $\mathbb{E}_{i}$, $\overbracket{\mathbb{E}_{ij}}=\Big\langle  (\overrightarrow{\eta_{ij}}, \overrightarrow{\xi_{ij}}); \overrightarrow{\tau^{\mathbb{T}}_{ij}}, \overrightarrow{\tau^{\mathbb{I}}_{ij}}, \overrightarrow{\tau^{\mathbb{F}}_{ij}} \Big\rangle  $ is aggregated into $\overbracket{\mathbb{E}_{i}}=\Big\langle  (\overrightarrow{\eta_{i}}, \overrightarrow{\xi_{i}}); \overrightarrow{\tau^{\mathbb{T}}_{i}}, \overrightarrow{\tau^{\mathbb{I}}_{i}}, \overrightarrow{\tau^{\mathbb{F}}_{i}}\Big\rangle  $.\\
{\bf Step4:} Compute the positive and negative ideal values of each alternatives. Here, the positive ideal value can be founded as,
$${\overbracket{\mathbb{E}^{+}}}= \Big\langle  \Big(\max_{1\leq i\leq n}(\overrightarrow{\eta_{ij}}), \min_{1\leq i\leq n}(\overrightarrow{\xi_{ij}})\Big); 1,1,0\Big\rangle$$ and
the negative ideal value can be founded as,
$${\overbracket{\mathbb{E}^{-}}}= \Big\langle  \Big(\min_{1\leq i\leq n}(\overrightarrow{\eta_{ij}}), \max_{1\leq i\leq n}(\overrightarrow{\xi_{ij}})\Big); 0,0,1\Big\rangle .$$
{\bf Step5:} Compute the Hamming distances between each alternative with two ideal values
$$\mathbb{D}_{i}^{+} = \mathbb{D}_{H} \Big(\overbracket{\mathbb{E}_i},{\overbracket{\mathbb{E}^{+}}}\,\Big);\,\, \mathbb{D}_{i}^{-} = \mathbb{D}_{H}\Big(\overbracket{\mathbb{E}_i},{\overbracket{\mathbb{E}^{-}}}\,\Big).$$
{\bf Step6:} Find the relative closeness values and find the ranking of alternatives
$$\mathbb{D}^{\ast}_{i}= \frac{\mathbb{D}_{i}^{-}}{\mathbb{D}_{i}^{+} + \mathbb{D}_{i}^{-}}$$
{\bf Step7:} Output yields for the optimal value is $\max \mathbb{D}^{\ast}_{i}$, and hence decision is to choose as the optimal solution to the problem.

\subsection{Selection of Successful Engineers}
Engineers are the inventors, designers, analyzers and builders of our modern age. They create the machines, structures systems etc. The limit of physics, the limit of modern age manufacturing technology, the limit imposed by current material properties, health and safety requirements, and costs - these are the things that engineers need to keep in mind when designing anything. A construction company is conscious to choose the most appropriate successful engineer to build a new construction work. After the first round of suggestions, a set of five engineers (alternatives) namely $\mathbb{E} =\{ \mathbb{E}_1, \mathbb{E}_2, \mathbb{E}_3,\mathbb{E}_4, \mathbb{E}_5 \}$ is suggested. An engineer is evaluated based on following attributes $\mathbb{A} = \{e_1 :$ {\text team work,} $e_2 $: {\text creativity}, $e_3$ : {\text analytical ability}, $e_4$ : {\text leadership}\} and their corresponding weights are $w = \{0.35, 0.27, 0.23, 0.15\}$. The decision-making informations are as follows: \\ 
{\bf Step 1:}   
According to the decision information  the following table is constructed.
\begin{center}
\begin{tabular}{||c c c c||} 
 \hline
$$ &$e_{1}$ &$e_{2}$ &$e_{3}$\\
 \hline \hline
$\mathbb{E}_{1} $  &  $\langle (0.85,0.5);0.88,0.8,0.8\rangle$ & $\langle (0.55,0.5);0.85,0.7,0.9\rangle$ &
$\langle (0.65,0.6);0.8,0.8,0.85\rangle$\\
$\mathbb{E}_{2} $  &  $\langle (0.7,0.65);0.75,0.95,0.75\rangle $ & $\langle (0.65,0.45);0.85,0.65,0.95\rangle $ &
$\langle (0.5,0.4);0.7,0.85,0.9\rangle $ \\
$\mathbb{E}_{3} $  &  $\langle (0.75,0.6);0.8,0.7,0.85\rangle $ & $\langle (0.75,0.7);0.9,0.85,0.8\rangle $ &
$\langle (0.65,0.5);0.75,0.9,0.8\rangle $  \\
$\mathbb{E}_{4} $  &  $\langle (0.5,0.4);0.7,0.85,0.9\rangle $ & $\langle (0.6,0.5);0.95,0.75,0.75\rangle $ &
$\langle (0.75,0.55);0.8,0.95,0.7\rangle $\\
$\mathbb{E}_{5} $  &  $\langle (0.65,0.55);0.9,0.75,0.85\rangle $ & $\langle (0.7,0.6);0.75,0.9,0.8\rangle $ & $\langle (0.7,0.65); 0.9,0.75,0.85\rangle $\\
[1ex] 
 \hline
\end{tabular}
\end{center}
\begin{center}
\begin{tabular}{||c c ||} 
 \hline
$$ &$e_{4}$\\
 \hline \hline
$\mathbb{E}_{1} $  &  $\langle (0.6,0.55);0.7,0.85,0.9\rangle$\\
$\mathbb{E}_{2} $  &  $\langle (0.8,0.65);0.85,0.95,0.65\rangle $ \\
$\mathbb{E}_{3} $  &  $\langle (0.65,0.5);0.8,0.7,0.95\rangle $ \\
$\mathbb{E}_{4} $  &  $\langle (0.75,0.7);0.75,0.85,0.85\rangle $ \\
$\mathbb{E}_{5} $  &  $\langle (0.7,0.5);0.85,0.9,0.7\rangle $ \\
[1ex] 
 \hline
\end{tabular}
\end{center}
{\bf Step 2:} Normalized decision matrix is constructed below.
\begin{center}
\begin{tabular}{||c c c c c||} 
\hline
$$ &$e_{1}$ &$e_{2}$ &$e_{3}$ &$e_{4}$\\
 \hline  \hline
$\mathbb{E}_{1} $  &  $\langle (1,0.4525);$ & $\langle (0.7333,0.6494);$ & $\langle (0.8667,0.8521);$ & $\langle 
(0.75,0.7202);$\\
$$  &  $0.88,0.8,0.8\rangle $ & $0.85,0.7,0.9\rangle $ & $0.8,0.8,0.85\rangle $ & $0.7,0.85,0.9\rangle $\\
$\mathbb{E}_{2} $  &  $\langle (0.8235,0.9286);$ & $\langle (0.8667,0.4451);$ & $\langle (0.6667,0.4923);$ & $\langle 
(1,0.7545);$ \\
$$  &  $0.75,0.95,0.75\rangle $ & $0.85,0.65,0.95\rangle $ & $0.7,0.85,0.9\rangle $ & $0.85,0.95,0.65\rangle $ \\
$\mathbb{E}_{3} $  &  $\langle (0.8824,0.7385);$ & $\langle (1,0.9333);$ & $\langle (0.8667,0.5917);$ & $\langle 
(0.8125,0.5495);$ \\
$$  &  $0.8,0.7,0.85\rangle $ & $0.9,0.85,0.8\rangle $ &$0.75,0.9,0.8\rangle $ & $0.8,0.7,0.95\rangle $ \\
$\mathbb{E}_{4} $  &  $\langle (0.5882,0.4923);$ & $\langle (0.8,0.5952);$ & $\langle (1,0.6205);$ & $\langle 
(0.9375,0.9333);$ \\
$$  &  $0.7,0.85,0.9 \rangle $ & $0.95,0.75,0.75\rangle $ & $0.8,0.95,0.7\rangle $ & $0.75,0.85,0.85\rangle $ \\
$\mathbb{E}_{5} $  &  $\langle (0.7647,0.716);$ & $\langle (0.9333,0.7347);$ &$\langle (0.9333,0.9286);$ & $\langle 
(0.875, 0.5102);$ \\
$$  &  $0.9,0.75,0.85\rangle $ & $0.75,0.9,0.8\rangle $ & $0.9,0.75,0.85\rangle $ & $0.85,0.9,0.7\rangle $ \\
[1ex] 
 \hline
\end{tabular}
\end{center}
{\bf Step 3:} Aggregate of the information with FNNWA operator of each alternative can be founded below.
\begin{center}
\begin{tabular}{||c c ||} 
\hline
$$ &$ FNNWA \,\,operator\,\,(\Lambda = 1)$\\
\hline\hline
$\overbracket{\mathbb{E}_{1}}$ & $\langle (0.8598,0.6377);0.8375,0.7863,0.8524 \rangle $\\
$\overbracket{\mathbb{E}_{2}}$ & $\langle (0.8256,0.6716);0.7924,0.8911,0.8160 \rangle $\\
$\overbracket{\mathbb{E}_{3}}$ & $\langle (0.9000,0.7290);0.8277,0.8068,0.8385 \rangle $\\
$\overbracket{\mathbb{E}_{4}}$ & $\langle (0.7925,0.6157);0.8441,0.8663,0.8017 \rangle $\\
$\overbracket{\mathbb{E}_{5}}$ & $\langle (0.8656,0.7391);0.8660,0.8299,0.8122 \rangle $\\
[1ex] 
 \hline
\end{tabular}
\end{center}
{\bf Step 4:} The positive and negative ideal values of the alternative can be founded below.
\begin{center}
\begin{tabular}{||c ||} 
\hline
${\overbracket{\mathbb{E}^{+}}}$ \\
\hline\hline
$\langle (0.9,0.6157),1,1,0 \rangle $ \\
[1ex] 
 \hline
\end{tabular}\,\,\,\,\,
\begin{tabular}{||c ||} 
\hline
${\overbracket{\mathbb{E}^{-}}}$ \\
\hline\hline
$\langle (0.7925,0.7391),0,0,1 \rangle $\\
[1ex] 
 \hline
\end{tabular}
\end{center}
{\bf Step 5:} The Hamming distance between each alternative with positive(respectively, negative) ideal values can be calculated below.
\begin{center}
\begin{tabular}{||c c c c c ||} 
\hline
$\mathbb{D}^{+}_{1}$ &$\mathbb{D}^{+}_{2}$ &$\mathbb{D}^{+}_{3}$ &$\mathbb{D}^{+}_{4}$ & $\mathbb{D}^{+}_{5}$\\
\hline\hline
$0.1954$ & $0.1746$ & $0.1776$ & $0.1759$ &  $0.1602$ \\
[1ex] 
\hline
\end{tabular} \,\,\,\,\,
\begin{tabular}{||c c c c c||} 
\hline
$\mathbb{D}^{-}_{1}$ &$\mathbb{D}^{-}_{2}$ &$\mathbb{D}^{-}_{3}$ &$\mathbb{D}^{-}_{4}$ &$\mathbb{D}^{-}_{5}$\\
\hline\hline
$0.1733$ & $0.1938$ & $0.1908$ & $0.1925$ &  $0.2082$ \\
[1ex] 
 \hline
\end{tabular}
\end{center}
{\bf Step 6:} We calculate the relative closeness value of each alternative below.
\begin{center}
\begin{tabular}{||c c c c c||} 
\hline
$\mathbb{D}^{*}_{1}$ &$\mathbb{D}^{*}_{2}$ &$\mathbb{D}^{*}_{3}$ &$\mathbb{D}^{*}_{4}$ &$\mathbb{D}^{*}_{5}$\\
 \hline  \hline
$0.4704$ & $0.5260$ & $0.5180$ & $0.5224$ &  $0.5651$ \\
[1ex] 
 \hline
\end{tabular}
\end{center}
{\bf Step 7:} Ranking of alternatives is done as $\mathbb{E}_{5} \geq \mathbb{E}_{2} \geq \mathbb{E}_{4} \geq \mathbb{E}_{3} \geq \mathbb{E}_{1}$.
Therefore, the most appropriate successful engineer $\mathbb{E}_{5}$ is the best option for the above construction job.

\subsection{Comparation for proposed and existing models}
We apply different methods to the above problem and we have the following tables.

\begin{center}
\begin{tabular}{||c c c c c ||} 
\hline
$\Lambda=1$ &$FNNWA$ &$FNNWG$ &$GFNNWA$ &$GFNNWG$\\
\hline\hline
$TOPSIS-\,Hamming$  &  $\mathbb{E}_{5}\geq \mathbb{E}_{2}\geq \mathbb{E}_{4} $ & $\mathbb{E}_{5}\geq \mathbb{E}_{3}\geq \mathbb{E}_{2} $ & $\mathbb{E}_{5}\geq \mathbb{E}_{2}\geq \mathbb{E}_{4} $ & $\mathbb{E}_{5}\geq \mathbb{E}_{3}\geq \mathbb{E}_{2} $\\
$distance\,(proposed)$ &  $\mathbb{E}_{3}\geq \mathbb{E}_{1}$ & $\mathbb{E}_{4}\geq \mathbb{E}_{1}$ & $\mathbb{E}_{3}\geq \mathbb{E}_{1}$ & $\mathbb{E}_{4}\geq \mathbb{E}_{1}$\\
$TOPSIS-\,Hamming$  &  $\mathbb{E}_{5}\geq \mathbb{E}_{1}\geq \mathbb{E}_{3} $ & $\mathbb{E}_{5}\geq \mathbb{E}_{1}\geq \mathbb{E}_{3} $ & $\mathbb{E}_{5}\geq \mathbb{E}_{1}\geq \mathbb{E}_{3} $ & $\mathbb{E}_{5}\geq \mathbb{E}_{1}\geq \mathbb{E}_{3} $\\
$distance \cite{palani1}$ &  $\mathbb{E}_{4}\geq \mathbb{E}_{2}$ & $\mathbb{E}_{4}\geq \mathbb{E}_{2}$ & $\mathbb{E}_{4}\geq \mathbb{E}_{2}$ & $\mathbb{E}_{4}\geq \mathbb{E}_{2}$\\
\hline
\end{tabular}
\end{center}
We get the following table by changing different values of $\Lambda$.
\begin{center}
\begin{tabular}{ccccccccccccccccccccc}
\midrule
\multicolumn{6}{c}{Relative closeness values} \\
\cmidrule(r){2-6}
$\Lambda$ & $\mathbb{D}^{*}_{1}$ &$\mathbb{D}^{*}_{2}$ &$\mathbb{D}^{*}_{3}$ &$\mathbb{D}^{*}_{4}$ &$\mathbb{D}^{*}_{5}$ &$Order$ \\
\midrule
$\Lambda=2$ & $ 0.4730	$ & $ 0.5311	$ & $ 0.5226	$ & $ 0.5316	$ & $ 0.5687$ & $\mathbb{E}_{5} \geq \mathbb{E}_{4} \geq \mathbb{E}_{2} \geq \mathbb{E}_{3} \geq \mathbb{E}_{1}$ &\\
$\Lambda=3$ & $ 0.4756	$ & $ 0.5362	$ & $ 0.5276	$ & $ 0.5406	$ & $ 0.5724$ & $\mathbb{E}_{5} \geq \mathbb{E}_{4} \geq \mathbb{E}_{2} \geq \mathbb{E}_{3} \geq \mathbb{E}_{1}$ &\\
$\Lambda=4$ & $ 0.478	$ & $ 0.5412	$ & $ 0.5328	$ & $ 0.5487	$ & $ 0.5759$ & $\mathbb{E}_{5} \geq \mathbb{E}_{4} \geq \mathbb{E}_{2} \geq \mathbb{E}_{3} \geq \mathbb{E}_{1}$ &\\
$\Lambda=5$ & $ 0.4804	$ & $ 0.5459	$ & $ 0.538	$ & $ 0.5559	$ & $ 0.5792$ & $\mathbb{E}_{5} \geq \mathbb{E}_{4} \geq \mathbb{E}_{2} \geq \mathbb{E}_{3} \geq \mathbb{E}_{1}$ &\\
$\Lambda=6$ & $ 0.4825	$ & $ 0.5502	$ & $ 0.5432	$ & $ 0.5622	$ & $ 0.5823$ & $\mathbb{E}_{5} \geq \mathbb{E}_{4} \geq \mathbb{E}_{2} \geq \mathbb{E}_{3} \geq \mathbb{E}_{1}$ &\\
$\Lambda=7$ & $ 0.4845	$ & $ 0.5541	$ & $ 0.5483	$ & $ 0.5677	$ & $ 0.5853$ & $\mathbb{E}_{5} \geq \mathbb{E}_{4} \geq \mathbb{E}_{2} \geq \mathbb{E}_{3} \geq \mathbb{E}_{1}$ &\\
$\Lambda=8$ & $ 0.4864	$ & $ 0.5576	$ & $ 0.5531	$ & $ 0.5726	$ & $ 0.588$ & $\mathbb{E}_{5} \geq \mathbb{E}_{4} \geq \mathbb{E}_{2} \geq \mathbb{E}_{3} \geq \mathbb{E}_{1}$ &\\
$\Lambda=9$ & $ 0.4881	$ & $ 0.5608	$ & $ 0.5576	$ & $ 0.5769	$ & $ 0.5906$ & $\mathbb{E}_{5} \geq \mathbb{E}_{4} \geq \mathbb{E}_{2} \geq \mathbb{E}_{3} \geq \mathbb{E}_{1}$ &\\
$\Lambda=10$ & $ 0.4897	$ & $ 0.5636	$ & $ 0.5617	$ & $ 0.5807	$ & $ 0.593$ & $\mathbb{E}_{5} \geq \mathbb{E}_{4} \geq \mathbb{E}_{2} \geq \mathbb{E}_{3} \geq \mathbb{E}_{1}$ &\\
$\Lambda=11$ & $ 0.4912	$ & $ 0.5661	$ & $ 0.5656	$ & $ 0.5841	$ & $ 0.5952$ & $\mathbb{E}_{5} \geq \mathbb{E}_{4} \geq \mathbb{E}_{2} \geq \mathbb{E}_{3} \geq \mathbb{E}_{1}$ &\\
$\Lambda=12$ & $ 0.4926	$ & $ 0.5683	$ & $ 0.5692	$ & $ 0.5872	$ & $ 0.5974$ & $\mathbb{E}_{5} \geq \mathbb{E}_{4} \geq \mathbb{E}_{2} \geq \mathbb{E}_{3} \geq \mathbb{E}_{1}$ &\\
$\Lambda=13$ & $ 0.4939	$ & $ 0.5704	$ & $ 0.5725	$ & $ 0.5901	$ & $ 0.5994$ & $\mathbb{E}_{5} \geq \mathbb{E}_{4} \geq \mathbb{E}_{3} \geq \mathbb{E}_{2} \geq \mathbb{E}_{1}$ &\\
$\Lambda=14$ & $ 0.4952	$ & $ 0.5722	$ & $ 0.5756	$ & $ 0.5928	$ & $ 0.6013$ & $\mathbb{E}_{5} \geq \mathbb{E}_{4} \geq \mathbb{E}_{3} \geq \mathbb{E}_{2} \geq \mathbb{E}_{1}$ &\\
$\Lambda=15$ & $ 0.4964	$ & $ 0.5739	$ & $ 0.5784	$ & $ 0.5952	$ & $ 0.603$ & $\mathbb{E}_{5} \geq \mathbb{E}_{4} \geq \mathbb{E}_{3} \geq \mathbb{E}_{2} \geq \mathbb{E}_{1}$ &\\
$\Lambda=16$ & $ 0.4975	$ & $ 0.5754	$ & $ 0.581	$ & $ 0.5975	$ & $ 0.6047$ & $\mathbb{E}_{5} \geq \mathbb{E}_{4} \geq \mathbb{E}_{3} \geq \mathbb{E}_{2} \geq \mathbb{E}_{1}$ &\\
$\Lambda=17$ & $ 0.4985	$ & $ 0.5768	$ & $ 0.5834	$ & $ 0.5997	$ & $ 0.6063$ & $\mathbb{E}_{5} \geq \mathbb{E}_{4} \geq \mathbb{E}_{3} \geq \mathbb{E}_{2} \geq \mathbb{E}_{1}$ &\\
$\Lambda=18$ & $ 0.4996	$ & $ 0.5781	$ & $ 0.5856	$ & $ 0.6017	$ & $ 0.6078$ & $\mathbb{E}_{5} \geq \mathbb{E}_{4} \geq \mathbb{E}_{3} \geq \mathbb{E}_{2} \geq \mathbb{E}_{1}$ &\\
$\Lambda=19$ & $ 0.5005	$ & $ 0.5793	$ & $ 0.5877	$ & $ 0.6036	$ & $ 0.6092$ & $\mathbb{E}_{5} \geq \mathbb{E}_{4} \geq \mathbb{E}_{3} \geq \mathbb{E}_{2} \geq \mathbb{E}_{1}$ &\\
$\Lambda=20$ & $ 0.5014	$ & $ 0.5804	$ & $ 0.5896	$ & $ 0.6054	$ & $ 0.6105$ & $\mathbb{E}_{5} \geq \mathbb{E}_{4} \geq \mathbb{E}_{3} \geq \mathbb{E}_{2} \geq \mathbb{E}_{1}$ &\\
$\Lambda=21$ & $ 0.5023	$ & $ 0.5815	$ & $ 0.5915	$ & $ 0.6071	$ & $ 0.6118$ & $\mathbb{E}_{5} \geq \mathbb{E}_{4} \geq \mathbb{E}_{3} \geq \mathbb{E}_{2} \geq \mathbb{E}_{1}$ &\\
$\Lambda=22$ & $ 0.5031	$ & $ 0.5824	$ & $ 0.5931	$ & $ 0.6087	$ & $ 0.613$ & $\mathbb{E}_{5} \geq \mathbb{E}_{4} \geq \mathbb{E}_{3} \geq \mathbb{E}_{2} \geq \mathbb{E}_{1}$ &\\
$\Lambda=23$ & $ 0.5039	$ & $ 0.5834	$ & $ 0.5947	$ & $ 0.6102	$ & $ 0.6141$ & $\mathbb{E}_{5} \geq \mathbb{E}_{4} \geq \mathbb{E}_{3} \geq \mathbb{E}_{2} \geq \mathbb{E}_{1}$ &\\
$\Lambda=24$ & $ 0.5047	$ & $ 0.5842	$ & $ 0.5962	$ & $ 0.6116	$ & $ 0.6151$ & $\mathbb{E}_{5} \geq \mathbb{E}_{4} \geq \mathbb{E}_{3} \geq \mathbb{E}_{2} \geq \mathbb{E}_{1}$ &\\
$\Lambda=25$ & $ 0.5055	$ & $ 0.585	$ & $ 0.5976	$ & $ 0.613	$ & $ 0.6161$ & $\mathbb{E}_{5} \geq \mathbb{E}_{4} \geq \mathbb{E}_{3} \geq \mathbb{E}_{2} \geq \mathbb{E}_{1}$ &\\
$\Lambda=26$ & $ 0.5062	$ & $ 0.5858	$ & $ 0.5989	$ & $ 0.6143	$ & $ 0.617$ & $\mathbb{E}_{5} \geq \mathbb{E}_{4} \geq \mathbb{E}_{3} \geq \mathbb{E}_{2} \geq \mathbb{E}_{1}$ &\\
$\Lambda=27$ & $ 0.5069	$ & $ 0.5865	$ & $ 0.6002	$ & $ 0.6156	$ & $ 0.6179$ & $\mathbb{E}_{5} \geq \mathbb{E}_{4} \geq \mathbb{E}_{3} \geq \mathbb{E}_{2} \geq \mathbb{E}_{1}$ &\\
$\Lambda=28$ & $ 0.5075	$ & $ 0.5871	$ & $ 0.6013	$ & $ 0.6168	$ & $ 0.6188$ & $\mathbb{E}_{5} \geq \mathbb{E}_{4} \geq \mathbb{E}_{3} \geq \mathbb{E}_{2} \geq \mathbb{E}_{1}$ &\\
$\Lambda=29$ & $ 0.5082	$ & $ 0.5878	$ & $ 0.6025	$ & $ 0.618	$ & $ 0.6195$ & $\mathbb{E}_{5} \geq \mathbb{E}_{4} \geq \mathbb{E}_{3} \geq \mathbb{E}_{2} \geq \mathbb{E}_{1}$ &\\
$\Lambda=30$ & $ 0.5088	$ & $ 0.5884	$ & $ 0.6035	$ & $ 0.6191	$ & $ 0.6203$ & $\mathbb{E}_{5} \geq \mathbb{E}_{4} \geq \mathbb{E}_{3} \geq \mathbb{E}_{2} \geq \mathbb{E}_{1}$ &\\
$\Lambda=31$ & $ 0.5094	$ & $ 0.589	$ & $ 0.6045	$ & $ 0.6201	$ & $ 0.621$ & $\mathbb{E}_{5} \geq \mathbb{E}_{4} \geq \mathbb{E}_{3} \geq \mathbb{E}_{2} \geq \mathbb{E}_{1}$ &\\
$\Lambda=32$ & $ 0.5099	$ & $ 0.5895	$ & $ 0.6055	$ & $ 0.6211	$ & $ 0.6217$ & $\mathbb{E}_{5} \geq \mathbb{E}_{4} \geq \mathbb{E}_{3} \geq \mathbb{E}_{2} \geq \mathbb{E}_{1}$ &\\
$\Lambda=33$ & $ 0.5105	$ & $ 0.59	$ & $ 0.6064	$ & $ 0.6221	$ & $ 0.6223$ & $\mathbb{E}_{5} \geq \mathbb{E}_{4} \geq \mathbb{E}_{3} \geq \mathbb{E}_{2} \geq \mathbb{E}_{1}$ &\\
$\Lambda=34$ & $ 0.5111	$ & $ 0.5905	$ & $ 0.6073	$ & $ 0.623	$ & $ 0.6229$ & $\mathbb{E}_{4} \geq \mathbb{E}_{5} \geq \mathbb{E}_{3} \geq \mathbb{E}_{2} \geq \mathbb{E}_{1}$ &\\
\midrule
\end{tabular}
\end{center}

\begin{center}
\includegraphics{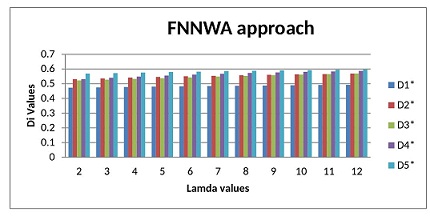}\\
{\bf Figure 1} Graphical representation consists of TOPSIS based on FNNWA approach.
\end{center}

\begin{center}
\includegraphics{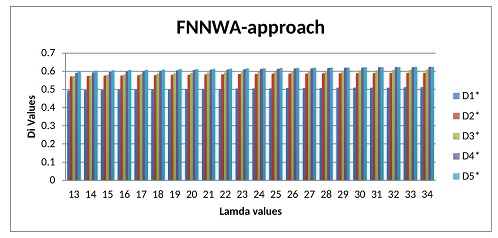}\\
{\bf Figure 2} Graphical representation consists of TOPSIS based on FNNWA approach.
\end{center}

From the above, we find that if $2 \leq \Lambda \leq 12$, then ranking of alternative is $\mathbb{E}_{5} \geq \mathbb{E}_{4} \geq \mathbb{E}_{2} \geq \mathbb{E}_{3} \geq \mathbb{E}_{1}$. If $13 \leq \Lambda \leq 33$, then ranking of alternative changes into $\mathbb{E}_{5} \geq \mathbb{E}_{4} \geq \mathbb{E}_{3} \geq \mathbb{E}_{2} \geq \mathbb{E}_{1}$.  If $\Lambda=34$, then ranking of alternative in a new order is $\mathbb{E}_{4} \geq \mathbb{E}_{5} \geq \mathbb{E}_{3} \geq \mathbb{E}_{2} \geq \mathbb{E}_{1}$. Hence the optimal alternative change from $\mathbb{E}_{5}$ into $\mathbb{E}_{4}$. Similarly, the alternative rankings are founded based on FNNWG, GFNNWA and FNNWG operator with the values of $\Lambda$.

\section{Conclusion:} \label{sec:7}
In this communication, we use euclidean and hamming distance measures for FNNSs, whose simplicity in the mathematical form is an extra boon to them. The supremacy of the euclidean and hamming distance measures is established via suitable numerical examples. In this example, the applicability of the euclidean and hamming distance measures is established. For the first time, we have proposed aggregation operations for FNNWA, FNNWG, GFNNWA, and GFNNWG operators and given some examples of how to develop these operators. The application of the FNNS MADM approach can help people make the correct decision out of available alternatives in indeterminate and inconsistent information environments. We have applied for the FNNWA, FNNWG, GFNNWA, and GFNNWG operators to the MADM problem on $\Lambda$. The distinct ranking of alternatives can be found with FNNWA, FNNWG, GFNNWA, and GFNNWG operators on $\Lambda$. In the finishing stage, the above analysis shows that the generalized values of $\Lambda$ have the greatest impact on the ranking of alternatives. The decision makers may set the values of $\Lambda$ according to the current situation for the best reasonable ranking, and then make appropriate decisions. Hence, the decision maker may determine the decision to arrive at the result based on the values of $\Lambda$. The euclidean and hamming distance measures of neutrosophic Fermatean sets have several real applications involving the analysis of data. The field of study is open, and the author is confident that the discussions in this article will be helpful to future researchers interested in this area of research.


\begin{thebibliography}{99}

\bibitem{Adeel}  A. Adeel, M. Akram and A.N.A. Koam, Group decision-making based on m-polar fuzzy linguistic TOPSIS method, Symmetry, 11(735), (2019), 1-20.

\bibitem{AkramArshad}  M. Akram and M. Arshad, A novel trapezoidal bipolar fuzzy TOPSIS method for group decision making, Group Decision and Negotiation, (2018), 1-20. 

\bibitem{Atanassov} K. Atanassov, Intuitionistic fuzzy sets, Fuzzy sets and Systems, 20(1), (1986), 87-96.

\bibitem{Ejegwa}  P.A. Ejegwa, Distance and similarity measures for Pythagorean fuzzy sets, Granular Computing, (2018), 1-17.

\bibitem{Hwang} C.L. Hwang, K. Yoon, Multiple attributes decision making methods and applications, Springer, Berlin Heidelberg, 1981.

\bibitem{Chiranjibe4} C. Jana and M. Pal, Multi criteria decision making process based on some single valued neutrosophic Dombi power aggregation operators, Soft Computing, 25(7), (2021), 5055-5072.

\bibitem{Jansi} R. Jansi, K. Mohana  and F. Smarandache, Correlation Measure for Pythagorean Neutrosophic Sets with $T$ and $F$ as Dependent Neutrosophic Components Neutrosophic Sets and Systems,30, (2019), 202-212.

\bibitem{Chiranjibe1} C. Jana, T. Senapati and M. Pal, Pythagorean fuzzy Dombi aggregation operators and its applications in multiple attribute decision making, International Journal of Intelligent Systems, 34(9), (2019), 2019-2038.

\bibitem{Chiranjibe7} C. Jana, G. Muhiuddin and M. Pal,  Multi criteria decision making approach based on SVTrN Dombi aggregation functions, Artificial Intelligence Review, 54(4), (2021), 3685-3723.

\bibitem{Chiranjibe8} C. Jana, M. Pal, F. Karaaslan  and J. Q. Wang,  Trapezoidal neutrosophic aggregation operators and their application to the multi-attribute decision making process, Scientia Iranica,  27(3), (2020), 1655-1673

\bibitem{Liang} D. Liang and Z. Xu, The new extension of TOPSIS method for multiple criteria decision 
making with hesitant Pythagorean fuzzy sets, Applied Soft Computing, 60, (2017), 167-179.    

\bibitem{palani1} M. Palanikumar, K. Arulmozhi, and C. Jana, Multiple attribute decision-making approach for Pythagorean neutrosophic normal interval-valued aggregation operators, Comp. Appl. Math. 41(90), (2022), 1-27.

\bibitem{PengDai}  X.D. Peng and J. Dai, Approaches to single-valued neutrosophic MADM based on MABAC, TOPSIS and new similarity measure with score function, Neural Computing and Applications, 29(10), (2018), 939-954.

\bibitem{Rahman1} K. Rahman, S. Abdullah,, M. Shakeel, MSA. Khan and M. Ullah, Interval valued 
Pythagorean fuzzy geometric aggregation operators and their application to group decision making problem, Cogent Mathematics, 4, (2017), 1-19.

\bibitem{Rahman2} K. Rahman, A. Ali, S. Abdullah and F. Amin,  Approaches to multi attribute group 
decision making based on induced interval valued Pythagorean fuzzy Einstein aggregation operator, New Mathematics and Natural Computation, 14(3), (2018), 343-361.

\bibitem{Shahzadi} G. Shahzadi, M. Akram and A. B. Saeid, An application of single-valued neutrosophic sets in medical diagnosis, Neutrosophic Sets and Systems, 18, (2017), 80-88.

\bibitem{Singh} P.K. Singh, Single-valued neutrosophic  context analysis at distinct multi-granulation. Comp. Appl. Math. 38, 80 (2019), 1-18. 

\bibitem{Senapati1} T. Senapati, R.R. Yager, Fermatean, fuzzy sets. J. Ambient Intell. Humaniz. Comput. 11, (2020), 663-674.

\bibitem{Smarandache1} F. Smarandache, A unifying field in logics, Neutrosophy neutrosophic probability, set and logic, American Research Press, Rehoboth, (1999).

\bibitem{Ali} K. Ullah, T. Mahmood, Z. Ali and N. Jan,  On some distance measures of complex Pythagorean fuzzy sets and their applications in pattern recognition, Complex and Intelligent Systems, (2019), 1-13.

\bibitem{Yager2} R.R. Yager, Pythagorean membership grades in multi criteria decision making, IEEE Trans. Fuzzy Systems, 22, (2014), 958-965.

\bibitem{Yang1} M.S Yang, C.H. Ko, On a class of fuzzy c-numbers clustering procedures for fuzzy data, Fuzzy Sets and Systems, 84, (1996), 49-60.

\bibitem{Zadeh} L. A. Zadeh,  Fuzzy sets, Information and control, 8(3), (1965), 338-353.

\bibitem{ZhangXu} X. Zhang and Z. Xu, Extension of TOPSIS to multiple criteria decision making with Pythagorean fuzzy sets, International Journal of Intelligent Systems, 29, (2014), 1061-1078.

\end{thebibliography}
\end{document}